\documentclass{amsart}
\usepackage{amscd}

\newtheorem{theorem}{Theorem}[section]
\newtheorem{lemma}[theorem]{Lemma}
\newtheorem{Cor}[theorem]{Corollary}
\theoremstyle{definition}
\newtheorem{definition}[theorem]{Definition}

\theoremstyle{remark}

\newcommand{\s}{\mathcal}

\newcommand{\sS}{\mathcal{S}}   
\newcommand{\scirc}{{\scriptstyle\circ}}
\newcommand{\Z}{\mathbb{Z}} 
\newcommand{\Q}{\mathbb{Q}}
\newcommand{\R}{\mathcal{R}}

\DeclareMathOperator{\im}{im}

\DeclareMathOperator{\Ext}{Ext} 
\DeclareMathOperator{\Hom}{Hom}
\DeclareMathOperator{\Sw}{Sw}
\DeclareMathOperator{\Ler}{Ler}
\DeclareMathOperator{\ASLer}{\skew{-18}\bar{\Ler}}
\DeclareMathOperator{\Ser}{LS}
\DeclareMathOperator{\Tor}{Tor}

\DeclareMathOperator{\fg}{fg}
\DeclareMathOperator{\Sd}{Sd}
\DeclareMathOperator{\inc}{inc}
\def\co{\colon\thinspace}

\begin{document}

\title{On Sikora's spectral sequences}

\author{Donald W. Barnes}

\address{1 Little Wonga Road, Cremorne, NSW 2090, Australia}
\email{donwb@iprimus.com.au}

\thanks{This work was done while the author was an Honorary Associate of the
School of Mathematics and Statistics, University of Sydney.}

\begin{abstract} In his study of Poincar\'e duality of a space $X$ acted on by a group
$G$, Sikora uses three spectral sequences which he calls the Leray, the Leray-Serre and the Swan
spectral sequences.  I show that, if $G$ is discrete and $X$ is  compact and Hausdorff, then the
Alexander-Spanier cohomology versions of these three spectral sequences are
isomorphic.  If also $X$ is HLC, then all versions of these spectral sequences are
isomorphic.  
\end{abstract}

\maketitle

\section{Introduction} \label{sec-Intro}
In his paper \cite{Sik}, Sikora establishes Poincar\'e duality properties for some spectral
sequences associated with the action of a Lie group $G$ on a space $X$.   
Let $EG \to BG$ be the universal $G$-bundle.  Let $X_G = X \times EG/\sim$ where
$\sim$ is the equivalence relation given by $(gx, e) \sim (x, g^{-1}e)$.  Let $R$ be a
commutative ring.  Sikora studied three spectral sequences associated with the bundle
$$\begin{CD}
X @>>>X_G = X \times EG/\sim\\
@.                      @V{\pi}VV \\
      @.                BG  \end{CD} $$
These spectral  sequences are
\begin{enumerate}
 \item The Leray spectral sequence of $\pi$ for the constant sheaf $\R = R \times X_G$
over $X_G$.  This is the composite functor spectral sequence for the functor
$\Gamma_{X_G} = \Gamma_{BG} \scirc \pi_*$, where $\Gamma$ denotes the global
section functor.
 \item The Leray-Serre spectral sequence of $\pi$ for the singular cohomology of $X_G$ with
coefficients in
$R$.
 \item The Swan\footnote{R. G. Swan informs me that this spectral sequence was
well-known before he adapted it in \cite{Swf} for application to fixed point theory.}
spectral sequence defined by Sikora to be the spectral sequence coming from the double
complex $\Sw^{pq} = \Hom(P_p, C^q)$, where $P_*$ is an $RG$-projective resolution of
$R$ and $C^\bullet$ is a cochain complex for the sheaf cohomology of the constant sheaf
$R \times X$ over $X$.
\end{enumerate}
Sikora defines and uses the Swan spectral sequence only in the case where $G$ is discrete.  It is
effectively limited to this case as the use of an $RG$-projective resolution ignores any topology
given on $G$.

Sikora has asked if the above spectral sequences are isomorphic.   For a discrete group $G$, I show
that, if the ring $R$ is noetherian and $X$ is a compact Hausdorff  space,  then the
Alexander-Spanier versions of the Leray, Leray-Serre and Swan spectral sequences are isomorphic. 
If also $X$ is homologically locally connected, then all versions of the spectral sequences are
isomorphic.   The Leray and Leray-Serre spectral sequences are defined for the cohomology of the
total spaces of arbitrary fibre bundles.  My analysis of the Leray-Serre spectral sequence requires  
that the bundle be definable with a discrete group of coordinate transformations,  so does not
answer Sikora's question for the spaces acted on by the circle group studied in his paper \cite{Sik}.

The Leray spectral sequence is defined for arbitrary sheaves over $X_G$.  So that the Leray-Serre
spectral sequence may be considered in comparable generality, in Section \ref{sec-coch}, I define
singular and Alexander-Spanier cohomology with coefficients in an arbitrary sheaf.  A version of this
is to be found in Bredon \cite{Bred2}.  The Swan spectral sequence can also be generalised.  This
is done using the concept of a $G$-sheaf, developed in Grothendieck \cite{Grot}.  I set out the
relevant theory in Section \ref{sec-Gsh}.

I am greatly indebted to R. G. Swan for the substantial contribution he has made to this
paper, pointing out errors in earlier versions, suggesting corrections, alternative
approaches and improvements.

\section{Summary of results}\label{sec-summary}
Let $R$ be a  noetherian ring and let $\s{A}$ be a sheaf of $R$-modules over the paracompact 
space $X$.  In Section \ref{sec-coch}, Alexander-Spanier and singular cohomology are
generalised to cohomology with coefficients in $\s{A}$.  A local  cochain is a cochain of some  
open subset $U \subset X$ with coefficients in the module $\s{A}(U)$ of sections over $U$. 
Taking germs of local cochains gives a sequence of sheaves over $X$.  The cohomology with
coefficients in $\s{A}$ is then defined to be the cohomology of the sequence of modules of   
global sections of these sheaves.  For each theory, we also have an fg version in which a local
cochain is required to have all its values in some finitely-generated submodule of $\s{A}(U)$.   
The following result is proved in Section \ref{sec-coch}.

\begin{theorem} \label{sum-coch} Let $R$ be a noetherian ring and let $\s{A}$ be a sheaf of
$R$-modules over the paracompact space $X$.  Then both the Alexander-Spanier and the fg
Alexander-Spanier cohomologies of $X$ with coefficients in $\s{A}$ are canonically isomorphic to
the sheaf cohomology of $\s{A}$.  Suppose also that $X$ is homologically locally connected. 
Then the singular and fg singular cohomologies of $X$ with coefficients in $\s{A}$ are also
canonically isomorphic to the sheaf cohomology of $\s{A}$. \end{theorem}

Let $G$ be a group acting on the space $X$ and compatibly on the sheaf $\s{A}$ over $X$.  We
call such a sheaf a $G$-sheaf.  The equivariant cohomology $H^n_G(X; \s{A})$ is defined to be
the cohomology of the sheaf $\s{A}_G = \s{A} \times_G EG$ over $X_G = X \times_G EG$.  The
main results on this, proved in Section \ref{sec-SwL}, are 

\begin{theorem}\label{th-eacyc} Let $G$ be a discrete group and let $\s{A}$ be a $G$-sheaf  
over the compact $G$-space $X$.  Suppose $\s{A}$ is $\Gamma^G_X$-acyclic.   Then $\s{A}_G$
is $\Gamma_{X_G}$-acyclic
\end{theorem} 
and
\begin{theorem}\label{th-equivar}  Let $G$ be a discrete group and let $\s{A}$ be a $G$-sheaf
over the compact $G$-space $X$.  Then
$$H^n_G(X; \s{A}) =  R^n\Gamma^G_X \s{A}.$$\end{theorem}

The purpose of this paper is to compare the Swan, Leray and Leray-Serre spectral
sequences associated with a space $X$ acted on by a group $G$.  The generalisations to
cohomology with coefficients in a $G$-sheaf were introduced to facilitate that comparison.  The
Leray and Swan spectral sequences are primarily defined in terms of sheaf cohomology.  The
spectral sequences from page 2 onward are independent of the choices of resolutions, so we can
use the Alexander-Spanier or fg Alexander-Spanier resolutions.  If we wish to use the singular or
fg singular resolutions, we need the space $X$ to be HLC.  In Section \ref{sec-SwL}, we prove

\begin{theorem} \label{th-SL} Let $G$ be a discrete group and let $\s{A}$ be a $G$-sheaf over
the compact Hausdorff $G$-space $X$.  Then the Leray and Swan spectral sequences for
Alexander-Spanier cohomology are isomorphic from page 2 onward.  If also $X$ is HLC, then the
Leray and Swan  spectral sequences for singular cohomology are isomorphic to them.
\end{theorem}

For the Leray-Serre spectral sequence, the situation is more complicated.  For each of
Alexander-Spanier and singular cohomology, we have two generalisations, one using the full version
of the cohomology, the other using the fg version.  For most spectral sequence
constructions, the effect of choices disappears at page 2, so we might expect  the full and fg
constructions to give the same spectral sequence from page 2 onwards, but I have not been able     
to prove this.  Consequently, in the following theorem, which follows from Theorem \ref{th-SerL}, I
have to specify the version of the sequence.  In the special case of constant coefficients $R$
considered in Sikora's paper \cite{Sik}, the two generalisations are the same.

\begin{theorem}  Let $G$ be a discrete group and let $\s{A}$ be a $G$-sheaf over the compact
Hausdorff space $X$. Then the Leray and the fg version of the Leray-Serre spectral sequences  for
Alexander-Spanier cohomology are isomorphic from page 2 onwards. If also $X$ is HLC, then also   
the Leray and the fg version of the Leray-Serre spectral sequences  for singular cohomology are
isomorphic to them.
\end{theorem}

\section{Preliminaries} \label{sec-prelims}

In this paper, all sheaves are sheaves of $R$-modules.  Most of the results require that the
base space of our sheaves be at least paracompact, (that is, Hausdorff and with every open
covering having a locally finite refinement)  so  we assume this throughout.  Where we must use a
space constructed from other data, the necessary paracompactness is ensured by the following
lemmas.

\begin{lemma} \label{lem-BGpara} Let $G$ be a discrete group.  Then its classifying space $BG$ is
paracompact and locally contractible.\end{lemma}

\begin{proof}  $BG$ is a CW complex.  By Lundell and Weingram \cite[Theorem II.4.2]{LW}, it is
paracompact.  By \cite[Theorem II.6.6]{LW}, it is locally contractible. \end{proof} 

\begin{lemma}\label{lem-Epara}  Let $E \to B$ be a fibre bundle.  Suppose that $B$
is paracompact and that the fibre $F$ is compact.  Then $E$ is paracompact. \end{lemma}

\begin{proof}  Let $\{W_\alpha\}$ be an open covering of $E$.  Every $W_\alpha$ is a union of
subsets of the form $U \times V$ where $U \subset B$ and $V \subset F$.  As we may replace
$\{W_\alpha\}$ by a refinement consisting of subsets of this form, we may suppose $W_\alpha =
U_\alpha \times V_\alpha$.  Take $b \in B$.  As $E_b$ is compact, it is covered by some finite
collection $W_1, \dots, W_n$ of the $W_\alpha$.  Put $U_b = \bigcap_{i=1}^n U_i$,
and $W_{b,i} = U_b \times V_i$.  Then $\{W_{b,i}\}$ is a refinement of the given covering of $E$.
Consider the covering $\{U_b\}$ of $B$.  As $B$ is paracompact, it has a locally finite refinement
$\{U'_\beta\}$.  Each $U'_\beta$ is contained in some $U_b$.  Put $W'_{\beta, i} = U'_\beta \times
V_i$.  Then $\{W'_{\beta, i}\}$ is a further refinement of the given covering of $E$.  We show that it
is locally finite.

Consider $e \in E$.  Then $b = \pi(e)$ has a neighbourhood $U$ which meets only finitely many
of the $U'_\beta$.  Thus $e$ has a neighbourhood which meets only finitely many of the
$W'_{\beta, i}$.   \end{proof}

\begin{Cor}\label{cor-XGpara} Suppose $G$ is a discrete group and that the $G$-space $X$ is
compact.  Then $X_G$ is paracompact. \end{Cor}

\begin{proof} By Lemma \ref{lem-BGpara}, $BG$ is paracompact, and, by assumption, the fibre
$X$ is compact.\end{proof}

To avoid trivial
complications, we assume of any presheaf $A$ that $A(\emptyset) = 0$.   Following Swan
\cite{Sw}, we denote the sheaf  of germs of the presheaf $A$ by $LA$. For a sheaf $\s{A}$
over the space
$X$, we denote the module of sections  over the open subset $U \subseteq X$ by
$\Gamma_U(\s{A})$ or by $\s{A}(U)$.  We denote by $P\s{A}$ the presheaf of these
modules of  sections.   A presheaf $A$ is the presheaf $P\s{A}$ of a sheaf if and only if it
satisfies the conditions
\begin{description}
\item[(S1) Assembly] For every covering $\{U_i\}$ of
$U$ by open subsets and elements $u_i \in A(U_i)$ whose restrictions agree on
intersections, there exists an element $u \in A(U)$ whose restriction to $U_i$ is $u_i$
for all $i$.
\item[(S2) Locally zero] The only locally zero element of $A(U)$ is $0$.
\end{description}

For any presheaf $A$, there is a homomorphism $A(U) \to (PLA)(U)$.  The kernel is the
set of locally zero elements of $A(U)$.  The homomorphism is an isomorphism if both (S1)
and (S2) are satisfied.

For much of the following, we shall need the assumption that the base space of our sheaves
is HLC, so I give the definition and some basic facts about HLC spaces here.

\begin{definition}\rm The space $X$ is said to be \textit{homologically locally connected
(HLC)} if for every $x \in X$ and every neighbourhood $U$ of $x$, there exists for each $n$ a
neighbourhood $V$ of  $x$, $V \subset U$ such that the map $i_n\co \tilde{H}_n(V) \to
\tilde{H}_n(U)$ induced by the inclusion $i\co V \to U$ is the zero map.  Here
$\tilde{H}_n(U)$ is the reduced singular homology group of $U$ with coefficients in
$\Z$. \end{definition}

Note that, if $i_n\co \tilde{H}_n(V) \to \tilde{H}_n(U)$ is the zero map and $V'
\subset V$ is another neighbourhood of $x$, then the map $i'_n\co \tilde{H}_n(V') \to
\tilde{H}_n(U)$ is also zero.  It follows that $V$ may be chosen such that $i_r = 0$ for
all $r \le n$.  Note also that the dimension $0$ part of the condition is that every
neighbourhood $U$ of $x \in X$ contains a neighbourhood $V$ such that every $x'
\in V$ can be joined to $x$ by a path in $U$.  This is equivalent to the path components of
$U$ being open and so, to  every neighbourhood $U$ of $x$ containing a path-connected
neighbourhood.  (Spanier \cite[p.\ 103 Exercise A1]{Sp}.)

Later, we will be working with spaces which, at least locally, are products.  So we
investigate the HLC property in relation to products.  

\begin{theorem} \label{th-HLCprod}  The product space $X \times Y$ is HLC if and only
if both $X$ and $Y$ are HLC. \end{theorem}

\begin{proof}  Since in the condition for a space to be HLC, we may always
replace a given neighbourhood by any smaller neighbourhood, we need only consider
neighbourhoods $W$ of $w = (x,y) \in X \times Y$ of the form $U \times V$ where
$U$ is a neighbourhood of $x$ in $X$ and $V$ is a neighbourhood of $y$ in $Y$.    As a
path from $(x,y)$ to $(x', y')$ in $U' \times V'$ is a pair of paths from $x$ to $x'$ in
$U'$ and from $y$ to $y'$ in $V'$, the condition in dimension $0$ holds for $X \times Y$
if and only if it holds for both $X$ and $Y$. 

Suppose $X$ and $Y$ are HLC.  Let $n >0$ and let $U$ and $V$ be neighbourhoods
of $x \in X$ and $y \in Y$.  Since $X$ is HLC, we can choose neighbourhoods $U_i$ of
$x$ such that $ U_2 \subset U_1 \subset U_0 = U$ and
such that the inclusions $U_i \to U_{i-1}$ induce the zero map $\tilde H_r(U_i) \to 
\tilde H_r(U_{i-1})$ for all $r \le n$.  We choose neighbourhoods $V_i$ of $y$
similarly.  Put $W_i = U_ i \times V_i$.  Since the K\"unneth short exact sequence is
natural, we have the commutative diagram with exact rows
$$\begin{CD}
\bigoplus\limits_{r+s=n} H_r(U_2) \otimes H_s(V_2) @>{\lambda_2}>> H_n(W_2) 
@>{\mu_2}>> \bigoplus\limits_{r+s = n-1} \Tor(H_r(U_2), H_s(V_2)) \\
   @VV{\theta_2}V   @VV{\phi_2}V   @VV{\psi_2}V \\
\bigoplus\limits_{r+s=n} H_r(U_1) \otimes H_s(V_1) @>{\lambda_1}>> H_n(W_1) 
@>{\mu_1}>> \bigoplus\limits_{r+s = n-1} \Tor(H_r(U_1), H_s(V_1)) \\
   @VV{\theta_1}V   @VV{\phi_1}V   @VV{\psi_1}V \\
\bigoplus\limits_{r+s=n} H_r(U_0) \otimes H_s(V_0) @>{\lambda_0}>> H_n(W_0) 
@>{\mu_0}>> \bigoplus\limits_{r+s = n-1} \Tor(H_r(U_0), H_s(V_0))   \end{CD} $$

Since $\theta_i$ is zero on any summand with either of $r, s$ non-zero and $r+s = n >
0$, we have $\theta_i = 0$.  Likewise, $\psi_i$ is zero on the summands with $r > 0$.  Since
$H_0(U_i)$ is free, $\Tor(H_0(U_i), H_s(V_i)) = 0$ and it follows that $\psi_i = 0$.  Thus
$\im\phi_2 \subseteq \ker\mu_1 = \im\lambda_1$.  As $\phi_1(\im\lambda_1) = 0$, we have
$\phi_1\phi_2 = 0$.

Conversely, suppose that $X \times Y$ is HLC.  By induction, we may suppose we have
$U'$ and $V'$ such that $\theta_r\co \tilde H_r(U') \to \tilde H_r(U)$ and $\phi_s\co
\tilde H_s(V') \to \tilde H_s(V)$ are zero for $r,s<n$.  From the naturality of the K\"unneth
Formulae, we have the commutative diagram
$$\begin{CD}
0 @>>> \bigoplus\limits_{r+s=n} H_r(U') \otimes H_s(V')  @>>> H_n(W')\\
@.      @VV{\lambda}V    @VV{\mu=0}V\\ 
0 @>>> \bigoplus\limits_{r+s=n}H_r(U) \otimes H_s(V)  @>>>H_n(W)
\end{CD}$$
where $\lambda = \bigoplus_{r+s=n} \theta_r \otimes \phi_s$.  Since $\mu=0$, we
have $\theta_0 \otimes \phi_n = 0$.   But $H_0 = \Z \oplus \tilde{H}_0$.   The
map $H_0(U') \to H_0(U)$ is the identity $\Z \to \Z$ and zero on $\tilde{H}_0(U')$.   It
follows that the map $\phi_n\co H_n(V') \to H_n(V)$ is zero.
\end{proof}

\section{Sheaves of cochains} \label{sec-coch}
In order to compare singular and sheaf cohomology, I consider sheaves of singular
cochains following Swan \cite[p.\  28]{Sw}.  Swan defines these with coefficients in a locally
constant sheaf.  As the  Leray spectral sequence is defined for arbitrary sheaves, I try to
generalise to coefficients in arbitrary sheaves.

A locally constant sheaf $\s{A}$ is the same as a system of local
coefficients\footnote{I use this only to motivate the definition of cohomology with
coefficients in a general sheaf.  The conditions ($X$ locally path connected and semi-locally
1-connected) under which it holds need not concern us.}.  A singular $n$-cochain $c$ with
coefficients in the sheaf $\s{A}$ assigns to each singular $n$-simplex $\sigma\co
\Delta^n \to X$, an element $c(\sigma) \in \s{A}_{\sigma(v_0)}$.  We denote the module of
all these $n$-cochains by $S^n(X; \s{A})$.  As $\s{A}$ is locally constant, any path in $X$
from $x_0$ to $x_1$ gives a map $\s{A}_{x_0} \to \s{A}_{x_1}$.  This enables us to shift
the values given by $c$ on the faces of an $(n+1)$-simplex to the first vertex and so define
$\delta c$.   We thus have the cochain complex $S^\bullet(X; \s{A})$ of singular cochains
with coefficients in $\s{A}$.  Following Swan, we obtain for each $n$, a presheaf $S^n$ by
setting $S^n(U) =  S^n(U; \s{A}|U)$.  We form its sheaf of germs, which we denote by
$\sS^n(X; \s{A})$.  First some facts about this.

 Each $c \in S^n(X,\s{A})$ gives a section $\phi c$ of $\s{S}^n$ by setting $\phi c(x)$ to be the
germ at $x$ of $c$.  This clearly defines a homomorphism $\phi\co S^n(X,\s{A}) \to \Gamma
\s{S}^n$.  Now suppose $\gamma \in \Gamma \s{S}^n$.  Each $\gamma(x)$ is the germ of
some cochain $c_x \in S^n(U_x; \s{A}|U_x)$ of some neighbourhood $U_x$ of $x$.  For
each $y \in U_x$, we have the germ $g(y)$ of $c_x$ at $y$, giving a section over $U_x$. 
The sections $\gamma$ and $g$ agree at $x$ and so on some neighbourhood $V_x$ of
$x$.  We thus have for each $x \in X$, a neighbourhood $V_x$ of $x$ and a cochain $c_x$
of $V_x$ such that the germ of $c_x$ at any $y \in V_x$ is $\gamma(y)$.  Our space $X$
is paracompact, so the covering $\{V_x\}$ has a locally finite refinement $\{U_\alpha \}$
with cochains $c_\alpha$ of $U_\alpha$ agreeing with $\gamma$ as above.  Suppose $x
\in U_\alpha$.  Since $c_x$ and $c_\alpha$ have the same germ at $x$, they agree on
some neighbourhood $V_{x, \alpha}$ of $x$.  Replacing $V_x$ by the intersection of the
$V_{x, \alpha}$ for the finitely many $\alpha$ with $x \in U_\alpha$, we may suppose
that $c_x(\sigma) = c_\alpha(\sigma)$ for every $\alpha$ with $x \in U_\alpha$.

We now construct a cochain $c \in S^n(X; \s{A})$.  Let $\sigma$ be any $n$-simplex of
$X$.  If there exist $\alpha$ such that $\sigma(\Delta^n) \subseteq U_\alpha$, choose any
one and set $c(\sigma) = c_\alpha(\sigma)$.  If no such $\alpha$ exists, set $c(\sigma) =
0$.  Then the germ at $x$ of $c$ is the germ at $x$ of $c_x$, that is, $\gamma(x)$.  Thus
$\phi$ is surjective.  Clearly $\ker \phi $ is the submodule $O^n(X,\s{A})$ of locally zero 
$n$-cochains.  The locally zero cochains clearly form a subcomplex of $S^\bullet$.   We
thus have $\Gamma \s{S}^\bullet (X; \s{A}) = S^\bullet(X; \s{A}) / O^\bullet(X,\s{A})$. 
The argument of  Swan \cite[p.\ 88, Proposition 6]{Sw} applies  to give
$H(O^\bullet) = 0$ and so $H(\Gamma\sS^\bullet) = H(S^\bullet)$.

\begin{definition}\rm  The sheaf $\s{A}$ over $X$ is called \textit{fine} if, for every
locally finite open covering $\{U_i\}$ of $X$, there exist endomorphisms $\theta_i$ of
$\s{A}$ such that the support $|\theta_i| \subseteq \bar U_i$ (closure of $U_i$) and $\sum_i
\theta_i = 1$.
\end{definition}

By Swan\cite[p.\ 84, Proposition 3]{Sw}, the sheaf $\sS^n$ is fine.  Note that,  by Swan\cite[p.\ 75, Corollary to
Proposition 5]{Sw}, fine sheaves
over paracompact spaces are $\Gamma$-acyclic. 
We have a sequence of sheaves
$$0 \to \s{A} \to \sS^0  \to \sS^1 \to \sS^2  \to \ldots.$$
If $X$ is HLC (homologically locally connected) this sequence is exact.   (I prove a
more general result, Lemma \ref{lem-HLCsim} below.)  Suppose $X$ is paracompact and
HLC.  Since the $\sS^n$ are $\Gamma$-acyclic, the sequence is a resolution which can be
used to calculate the sheaf cohomology of $\s{A}$.   It follows that the sheaf cohomology of
$\s{A}$  coincides with the singular cohomology of $X$ with coefficients in $\s{A}$.  

Now  to generalise to sheaves which are not locally constant.  The definition of
$\s{S}^n(X; \s{A})$ given above does not lend itself to generalisation to sheaves which are    
not locally constant, so I give an equivalent formulation, beginning with a simple but useful
observation.

\begin{lemma}\label{constsec}  Let $A \times U$ be a constant sheaf over $U$ and let $V$ be
a connected subset of $U$. Then every section of $A \times U$  is constant on $V$.  More
generally, let
$\s{A}$ be a locally constant sheaf over the locally path-connected space $X$.  Then any
neighbourhood $U$ of $x \in X$ contains a neighbourhood $V$ such that every section of
$\s{A}$ over $V$ is constant.\end{lemma}

\begin{proof}  A section over $U$ is a continuous function into the discrete space $A$,
so is constant on connected subsets.  Any point $x \in X$ has a neighbourhood $U$ such that
$\s{A}|U$ is constant.  But $U$ contains a path-connected neighbourhood $V$ of $x$.
\end{proof}

In constructing germs of cochain with values in a locally constant sheaf, we can restrict our
attention to neighbourhoods $U_x$ of $x$ such that $\s{A}|U_x$ is constant.  Then a cochain    
of $U_x$ is a function which assigns to each simplex of $U_x$ a constant section over $U_x$.   
The tracking from one stalk to another required for the coboundary is now taken care of by the
restriction maps.  A general sheaf need not have such constant sections, so to generalise, we  
must allow arbitrary sections over $U_x$.  If the neighbourhood $U_x$ is not connected, the  
locally constant sheaf may have non-constant sections over $U_x$.  However, if the space is    
HLC in dimension $0$, then by Lemma \ref{constsec}, any neighbourhood $U_x$ of $x$
contains a neighbourhood $V_x$ over which, all sections are constant.  It follows that germs of
cochains of $U_x$ with values in $\s{A}(U_x)$ are germs of cochains of $V_x$ with values in   
 the module of constant sections over $V_x$.   Thus allowing arbitrary sections does not change
our sheaves of germs.

 For any sheaf $\s{A}$ and open set $U \subseteq X$, we have the module
$\s{A}(U)$ of sections over $U$.  We can form the cochain module $S^n(U; \s{A}(U))$.  If $V
\subseteq U$ is open, we have a restriction map $\rho\co S^n(U; \s{A}(U)) \to S^n(V; \s{A}(V))$
given by the restriction of the cochains to $V$ compounded with the map $\s{A}(U) \to
\s{A}(V)$.  This clearly defines a presheaf and we can take its sheaf $\sS^n(X; \s{A})$ of
germs.   We have a coboundary operator $\delta \co S^n(U; \s{A}(U)) \to S^{n+1}(U; \s{A}(U))$ 
and so obtain a coboundary operator $\delta \co \sS^n (X; \s{A}) \to \sS^{n+1} (X; \s{A})$.  We
can then define the singular cohomology of $X$ with coefficients in $\s{A}$ to be the
cohomology of the cochain complex
$$0 \to \Gamma \sS^0(X;\s{A}) \to \Gamma \sS^1(X;\s{A}) \to \Gamma \sS^2(X;\s{A}) 
\to \ldots.$$
If $\s{A}$ is locally constant, this is the definition given in the discussion above.
The sequence of sheaves $\s{S}^\bullet(X; \s{A})$ has an augmentation $\epsilon\co \s{A}
\to \s{S}^0(X; \s{A})$.  Suppose $a \in \s{A}_x$.  Then $a$ is the germ of some section
$s$ over some neighbourhood $U$ of $x$.  We can identify $s$ with the function
$\sigma\co U \to \s{A}(U)$ with $\sigma(u) = s$ for all $u \in U$.  Define $\epsilon(a)$ to
be the germ of $\sigma$.  Then $\epsilon$ is clearly injective.  Note that as 
$\sigma(\partial c) = 0$ for any $1$-chain of $U$, the germ of $\sigma$ is a
$0$-cocycle.

Bredon, in his second edition \cite[section III.1]{Bred2}, gives a different definition of
singular cohomology with coefficients in a sheaf.  He defines the singular cohomology of
$X$ with coefficients in $\s{A}$ to be the cohomology of the cochain complex $\Gamma
(\s{S}^n(X; \Z) \otimes_\Z \s{A})$ that is, of $ \Gamma (\s{S}^n(X; \s{Z}) \otimes_\Z
\s{A})$, where $\s{Z}$ is the constant sheaf $\Z \times X$.  As we will be working with
sheaves of $R$-modules, we consider the complex $\Gamma (\s{S}^n(X; R) \otimes_R
\s{A}) = \Gamma (\s{S}^n(X; \s{R}) \otimes_R \s{A})$ where $\s{R} = R \times X$.

\begin{lemma}\label{lem-fg} Let $\s{A}$ be a sheaf on $X$.  Then there is a
homomorphism 
$$i\co  \s{S}^n(X; \s{R}) \otimes_R \s{A}\to \s{S}^n(X; \s{A}).$$
The image of $i$ is the set of germs of cochains $c \in S^n(U; A_c)$, where $A_c$ is a
finitely generated submodule of $\s{A}(U)$.  If $R$ is noetherian, then $i$ is injective.
\end{lemma}

\begin{proof} An element of $\s{S}^n(X; \s{R}) \otimes_R \s{A}$ is the germ  of
an element of $S^n(U; R) \otimes_R \s{A}(U)$ for some open set $U \subseteq X$, thus the
germ of $c = \sum_{j=1}^k c_j \otimes a_j$ where $c_j \in S^n(U; R)$ and $a_j \in
\s{A}(U)$.  Let $A_c$ be the submodule of $\s{A}(U)$ generated by $ a_1, \ldots, a_k$. 
Then for every $n$-simplex $\sigma$ of $U$, $c(\sigma) =  \sum_{j=1}^k c_j(\sigma) a_j
\in A_c$.  Conversely, let $A_c$ be generated by $a_1, \ldots, a_k$ and let $c
\in S^n(U; \s{A}(U))$ have all its values in $A_c$.  Then, for each $n$-simplex $\sigma$,
$c(\sigma) = \sum_{j=1}^k c_j(\sigma) a_j$, where $c_j(\sigma) \in R$.  Choosing such a
representation for each $\sigma$ gives $c_j \in S^n(U; R)$ and the germ of $c$ is the
image under $i$ of the germ of $\sum_{j=1}^k c_j \otimes a_j$.

That $i$ is injective is essentially Cartan and Eilenberg \cite[Chapter II, Ex. 2]{CE}
\end{proof}

In view of this finite-generation property, I always assume from here on, that $R$ is
noetherian and identify $\s{S}^n(X; \s{R}) \otimes_R\s{A}$ with its image under $i$ which
I denote by
$\s{S}_{\fg}^n(X; \s{A})$.   Note that $S^n(U; R) \otimes_R \s{A}(U) \neq   S^n(U;
\s{A}(U))$ in general.   The values of a cochain in $S^n(U; \s{A}(U))$ need not lie in
any finitely generated submodule of $ \s{A}(U)$.  However, the inclusion
$\s{S}_{\fg}^\bullet(X; \s{A}) \to \s{S}^\bullet(X; \s{A})$ induces an isomorphism of the
cohomology modules, at least for paracompact HLC spaces for which, as we shall see in
Theorem \ref{th-canon} below, both coincide with the sheaf cohomology.  But the argument  
given above relating $ H(\Gamma \s{S}^\bullet(X; A))$ to $H(S^\bullet(X; A))$ for constant
coefficients $A$  fails when applied to $\Gamma\s{S}_{\fg}(X; A)$ and $S^\bullet(X; R)
\otimes_R A$.  Assembling local cochains, each with values in a finitely generated subgroup
of  $A$, may give a global cochain whose values do not lie in any finitely generated subgroup.  
Indeed, if $X$ has infinitely many path components and if $A$ is infinitely generated, then
$H^0(S^\bullet(X; R) \otimes A) \ne H^0(S^\bullet(X; A))$.  Of  course, if $A$ is finitely
generated, then $S^\bullet(X; R) \otimes A = S^\bullet(X; A)$.  If the space $X$ is
compact, then the assembly puts together finitely many local cochains, so the result is
again a cochain with values in a finitely generated subgroup and the argument works.

To simplify notation, I write $\otimes$ for $\otimes_R$ and $\Hom$ for $\Hom_R$ from
here on.  If operations over $\Z$ are required as in the proof of the next lemma, it will be
indicated.

\begin{lemma} \label{lem-HLCsim} Suppose $X$ is HLC.  Then the sequence
$$0 \to \s{A} \to \s{S}^0(X;\s{A}) \to  \s{S}^1(X;\s{A}) \to \s{S}^2(X;\s{A}) \to \ldots$$
is exact.  \end{lemma}

\begin{proof}  We have already seen that the sequence is exact at $\s{A}$.
We now show that the sequence is exact at $\s{S}^0(X; \s{A})$.   Let $\sigma \in
\s{S}^0_x$ and suppose that $\delta \sigma = 0$.  Then $\sigma$ is the germ of some
not necessarily continuous function $s\co U \to \s{A}(U)$ and $s(u_1) = s(u_2)$ for any
$u_1, u_2$ which are joined by a path in $U$.  Now let $V \subset U$ be a neighbourhood
of $x$ as in the definition of HLC.  Then every $v \in V$ can be joined to $x$ by a path in
$U$, so $s(v) = s(x)$ for all $v \in V$.  Thus the restriction $\rho s \in S^0(V; \s{A}(V))$
is constant and its germ
$\sigma$ is the germ of a section over $V$ of $\s{A}$.  Thus the sequence is exact at
$\s{S}^0$. 

Now consider $n > 0$.  Let $\sigma \in \s{S}^n_x$ and suppose $\delta \sigma = 0$. 
Then $\sigma$ is the germ of some $s \in S^n(U; \s{A}(U))$, $\delta s = 0$, where $U$ is
a neighbourhood of $x$.  There exist neighbourhoods $U_2 \subset U_1 \subset U_0 = U$
of $x$ such that the maps $i^k_n, i^k_{n-1}$ in reduced homology induced by the
inclusions $i^k\co U_k \to U_{k-1}$ are all $0$.  In the commutative diagram 
$$\begin{CD}
\Ext_\Z(H_{n-1}(U), A) @>>> H^n(U; A) @>>> \Hom_\Z(H_n(U), A) \\
@VVV  @VVV  @VVV\\
\Ext_\Z(H_{n-1}(U_1), A) @>>> H^n(U_1; A) @>>> \Hom_\Z(H_n(U_1), A) \\
@VVV  @VVV  @VVV\\
\Ext_\Z(H_{n-1}(U_2), A) @>>> H^n(U_2; A) @>>> \Hom_\Z(H_n(U_2), A) \\
\end{CD}$$
the rows are exact and the outside vertical arrows are all $0$.  It follows that the composite
$H^n(U; A) \to H^n(U_2; A)$ of the central arrows is $0$.  Taking $A=\s{A}(U)$ and
combining with the restriction $\s{A}(U) \to  \s{A}(U_2)$, we get that the restriction
$\rho s \in S^n(U_2; \s{A}(U_2))$ is a coboundary.  Thus there exists $t \in S^{n-1}(U_2;
\s{A}(U_2))$ with $\delta t = \rho s$.  Passing to germs, we obtain $\tau \in
\s{S}^{n-1}_x$ with $\delta \tau = \sigma$. 
 \end{proof} 

\begin{lemma} \label{lem-exfg} Suppose $X$ is HLC and that $R$ is noetherian.  Then the
sequence
$$0 \to \s{A} \to \s{S}_{\fg}^0(X;\s{A}) \to  \s{S}_{\fg}^1(X;\s{A}) \to
\s{S}_{\fg}^2(X;\s{A}) \to \ldots$$ is exact.  \end{lemma}

\begin{proof}   For each $a \in \s{A}_x$, the germ of a section $s \in \s{A}(U)$, we have the
\mbox{$0$-cochain} $c \in S^0(U; \s{A}(U))$ given by $c(u) = s$ for all $u \in U$.  This cochain
has values in the \mbox{$1$-generator} submodule $\langle s \rangle$ of $\s{A}(U)$.  Thus
$\epsilon(a) \in \s{S}_{\fg}^0(X; \s{A})$.  As $\epsilon$ is injective, the sequence is exact
at $\s{A}$.  By Lemma \ref{lem-HLCsim}, it is exact at $\s{S}_{\fg}^0$.

Now suppose $\sigma \in (\s{S}_{\fg}^n)_x$, $n > 0$ and that $\delta \sigma = 0$.  Then
$\sigma$ is the germ of some $s \in S^n(U; A)$, $\delta s = 0$, where $A$ is a finitely
generated submodule of $\s{A}(U)$.  As in the proof of Lemma \ref{lem-HLCsim}, we
obtain $U_2 \subset U$  and $t \in S^{n-1}(U_2; A)$ with $\delta t = s | U_2$.  Passing to
germs, we have that the germ $\tau$ of $t$ is in $\s{S}^{n-1}_{\fg}$ and $\delta \tau =
\sigma$.
\end{proof}

Similar constructions are possible with Alexander-Spanier cochains.  I put bars over
symbols to distinguish Alexander-Spanier (A-S) objects from the corresponding singular
objects.  For any abelian group $A$, an A-S $n$-cochain of $X$ is a function $c\co X^{n+1}
\to A$.  Denote the group of A-S $n$-cochains by $\bar S^n(X; A)$.  For a sheaf $\s{A}$,
we have the presheaf $\bar S^n(U; \s{A}(U))$ and form its sheaf of germs $\bar{\s{S}}(X;
A)$.  The A-S cohomology of $X$ with coefficients in $\s{A}$ is the cohomology of the
cochain complex $\Gamma_X \bar{\s{S}}^\bullet(X; \s{A})$.  For $\s{A}$ locally
constant, this is the standard definition of Alexander-Spanier cohomology.  We put
$\bar S_{\fg}^n(U; \s{A}) = \bar S^n(U; R) \otimes \s{A}(U)$ (assuming as always, that
$R$ is noetherian) and define
$\bar{\s{S}}_{\fg}^n(X; \s{A})$ to be the sheaf of germs of this presheaf.  Thus
$\bar{\s{S}}_{\fg}^n(X; \s{A}) = \bar{\s{S}}^n(X; \s{R}) \otimes \s{A}$ and  
$H^\bullet(\Gamma_X \bar{\s{S}}_{\fg}^\bullet(X; \s{A}))$ is Bredon's definition of
Alexander-Spanier cohomology with coefficients in $\s{A}$, given in \cite[section
III.2]{Bred2}.

The sheaves $\bar{\s{S}}^n(X; \s{A})$ are closely, but awkwardly, related to the sheaves
of the canonical flabby resolution of $\s{A}$.  The elements of that can be regarded as
germs of functions of $n+1$ variables, $f(x_0, \dots, x_n) \in \s{A}_{x_n}$, with $x_i$
in a subset depending on $x_0, \dots x_{i-1}$ and with $f(x_0, \dots, x_n) = 0$ if $x_1
= x_0$.  This is discussed in Godement \cite[Remarque 4.3.2]{God}.

\begin{lemma}  The sequences
$$0\to \s{A} \to \bar{\s{S}}^0(X;\s{A}) \to \bar{\s{S}}^1(X;\s{A}) \to \bar{\s{S}}^2
(X;\s{A}) \to \ldots$$ 
and
$$0\to \s{A} \to \bar{\s{S}}_{\fg}^0(X;\s{A}) \to \bar{\s{S}}_{\fg}^1(X;\s{A}) \to 
\bar{\s{S}}_{\fg}^2 (X;\s{A}) \to \ldots$$ 
are exact.  \end{lemma}

\begin{proof}  For any $R$-module $A$, the sequence
$$0 \to A \to \bar S^0(X; A) \to \bar S^1(X; A) \to \bar S^2(X; A) \to \ldots$$ is exact.  It
follows that the sheaf sequence 
$$0\to \s{A} \to \bar{\s{S}}^0(X;\s{A}) \to \bar{\s{S}}^1(X;\s{A}) \to \bar{\s{S}}^2
(X;\s{A}) \to \ldots$$
is also exact.  A cocycle $s \in \bar{S}_{\fg}^n(U;  \s{A}(U))$ with values in the finitely
generated submodule $A \subseteq \s{A}(U)$ is, when restricted to some smaller
neighbourhood $U'$, the coboundary of some $t \in S^{n-1}(U'; A)$.  Passing to germs, we
get that the second sequence is also exact.
\end{proof}

\begin{lemma}\label{lem-fine} The sheaves $\s{S}^n(X; \s{A})$, $\bar{\s{S}}^n(X; \s{A})$,
$\s{S}_{\fg}^n(X; \s{A})$ and $\bar{\s{S}}_{\fg}^n(X; \s{A})$ are fine. 
\end{lemma}

\begin{proof} Following Swan \cite{Sw}, let $\Delta$ be a space and let $v \in \Delta$ be
a selected point.  For singular cohomology, we use $\Delta = \Delta^n$, the standard
$n$-simplex, $v = v_0$ its first vertex, while for Alexander-Spanier cohomology, we
take $\Delta = \{v_0, \ldots, v_n\}$ and $v = v_0$.  For the open set $U \subseteq X$,
let $F(U)$ be the set of maps $\Delta \to U$ and let $T(U)$ be the $R$-module of all
functions $F(U) \to \s{A}(U)$.  Let $\s{S}$ be the sheaf of germs of this presheaf.  We
prove that $\s{S}$ is fine.

Let $\{U_i\}$ be a locally finite open covering of $X$.  Then each point $x \in X$ is in
finitely many of the $U_i$.  For each $x$, we choose one of these.  Denote the chosen
one by $U(x)$.  Thus $x \in U(x) \in \{U_i\}$.  To define the $\theta_i$, we must define
$\theta_i(s) \in \s{S}_x$ for each $s \in \s{S}_x$.  But $s$ is the germ at $x$ of some
function $c\co F(U) \to \s{A}(U)$ for some open set $U$.  We define $\theta_i(c)$ and take
$\theta_i(s)$ to be the germ at $x$ of $\theta_i(c)$.  To do this, for $\sigma \in F(U)$, we
put
$$\theta_i(c)(\sigma) = \begin{cases} c(\sigma) &\text{if $U(\sigma(v)) = U_i$}\\
                                     0 &\text{otherwise.} \end{cases}$$
Then $\theta_i$ is clearly an endomorphism.  If $x \not\in \bar{U}_i$, then the
neighbourhood $U$ of $x$ used above may be chosen disjoint from $U_i$ ensuring that
$\theta_i(c)(\sigma) = 0$.  Thus $|\theta_i| \subseteq \bar{U}_i$.   Further, $U(\sigma(v))
= U_i$ for exactly one $i$, so $\sum_i \theta_i(c)(\sigma) = c(\sigma)$.

Thus $\s{S}^n(X; \s{A})$ and $\bar{\s{S}}^n(X; \s{A})$ are fine.  
By the above, $\sS^n(X; \s{R})$ and $\bar{\sS}^n(X; \s{R})$ are fine.  By Spanier
\cite[p.\ 331, Theorem 3]{Sp}, $\sS^n(X; \s{R}) \otimes \s{A}$ and $\bar{\sS}^n(X; \s{R})
\otimes \s{A}$ are also fine. \end{proof}

If $X$ is paracompact, then a fine sheaf over $X$ is $\Gamma_X$-acyclic.  Thus $\bar{\s
S}^\bullet(X; \s{A})$ and $\bar{\s{S}}_{\fg}^\bullet(X; \s{A})$ are resolutions of $\s{A}$ by
$\Gamma_X$-acyclic sheaves.  If also $X$ is HLC, then so are $\s{S}^\bullet(X; \s{A})$
and $\s{S}_{\fg}^\bullet(X; \s{A})$.
We now have:
\begin{theorem} \label{th-canon} Let $X$ be a paracompact space. Suppose that 
$R$ is noetherian.  Then   $H^\bullet(\Gamma_X \bar{\s{S}}^\bullet(X; \s{A}))$, and
$H^\bullet(\Gamma_X \bar{\s{S}}_{\fg}^\bullet(X; \s{A}))$ are canonically isomorphic to
the sheaf cohomology $H^\bullet(X; \s{A})$.  If $X$ is HLC, then
$H^\bullet \Gamma_X (\s{S}^\bullet(X; \s{A}))$ and $H^\bullet(\Gamma_X
\s{S}_{\fg}^\bullet(X; \s{A}))$ are also canonically isomorphic to $H^\bullet(X; \s{A})$. 
\end{theorem}

We now compare $\s{S}_{\fg}^\bullet(X; \s{A})$ and $\s{S}^\bullet(X; \s{A})$.  As we saw
above, each element $c \in S^n(U; R) \otimes \s{A}(U)$ gives, for each singular
$n$-simplex of $U$, an element of $\s{A}(U)$ and so may be regarded as an element of
$S^n(U; \s{A}(U))$.  This inclusion passes to germs, and we have an inclusion $i\co
\s{S}_{\fg}^n(X; \s{A})  \to \s{S}^n(X; \s{A})$.  Even for locally constant sheaves, this
inclusion is not an isomorphism.

An Alexander-Spanier cochain $c$ can be regarded as a singular cochain by defining for
the singular simplex $\sigma$, $c(\sigma) = c(\sigma(v_0), \ldots, \sigma(v_n))$.  We
thus have a diagram of inclusions
$$\begin{CD}
\bar{\s{S}}_{\fg}^\bullet(X; \s{A})    @>>>   \bar{\s{S}}^\bullet(X; \s{A})\\
@VVV                                                    @VVV\\
\s{S}_{\fg}^\bullet(X; \s{A})             @>>>   \s{S}^\bullet(X; \s{A})\\
\end{CD}. $$
Taking global sections  gives a diagram of inclusions which, by Theorem
\ref{th-canon}, induce isomorphisms in cohomology if $X$ is paracompact
and HLC.

\begin{lemma} \label{lem-exfloc} Let $\s{A} \xrightarrow{\alpha} \s{B}
\xrightarrow{\beta} \s{C}$ be an  exact sequence of locally constant sheaves over a locally
path connected space $X$.  Then the sequences  
$$ \s{S}^n( X; \s{A}) \xrightarrow{\alpha} \s{S}^n(X; \s{B}) \xrightarrow{\beta}
\s{S}^n(X;
\s{C}) \text{ and } \bar{\s{S}}^n( X; \s{A}) \xrightarrow{\alpha} \bar{\s{S}}^n(X; \s{B})
\xrightarrow{\beta} \bar{\s{S}}^n(X; \s{C})$$
 are exact.
\end{lemma}

\begin{proof}  Let $g \in \s{S}^n(X; \s{B})$ and suppose $\beta(g) = 0$.  Then $g$ is the
germ at $x \in X$ of some cochain $b \in S^n(U; \s{B}(U))$ for some neighbourhood $U$ of
$x$, and the cochain $\beta b \in  S^n(U; \s{C}(U))$ has germ $0$.  $U$ may be
chosen such that the restrictions of $\s{A}, \s{B}$ and $\s{C}$ to $U$ are constant.  Thus
for some neighbourhood $V \subseteq U$ of $x$, $\beta b| V = 0$.  Since $X$ is locally path
connected, there exists a path connected neighbourhood $W \subseteq V$.  By Lemma
\ref{constsec}, all sections of
$\s{A}, \s{B}$ and $\s{C}$ over $W$ are constant sections.   For
any $n$-simplex $\sigma$ of $W$, $b(\sigma)$ is a constant section over $W$ and we have
$\beta b(\sigma) = 0$. Since the sequence of constant sections over $W$ is exact, for each
$\sigma$, we can choose a constant section $a(\sigma)$ of $\s{A}$ over
$W$ such that $\alpha a(\sigma) = b(\sigma)$.  This defines a cochain $a \in S^n(W;
\s{A}(W))$ such that $\alpha a = b$.  The germ of $a$ maps to $g$.  Thus 
$\s{S}^n( X; \s{A}) \xrightarrow{\alpha} \s{S}^n(X; \s{B}) \xrightarrow{\beta} \s{S}^n(X;
\s{C})$  is exact.  Similarly, $\bar{\s{S}}^n( X; \s{A}) 
\xrightarrow{\alpha} \bar{\s{S}}^n(X; \s{B}) \xrightarrow{\beta} \bar{\s{S}}^n(X; \s{C})$
is exact.
\end{proof}

Note that the use of Alexander-Spanier cochains here does not obviate the need to
assume that $X$ is locally path connected (HLC in dimension $0$).

\begin{lemma} \label{lem-Brexact} Suppose that $R$ is noetherian.  Then $\s{S}_{\fg}^n(X;
\_)$ is an exact functor on the category of sheaves over $X$. \end{lemma}
 \begin{proof} This is essentially Cartan and Eilenberg \cite[Chapter VI, Ex. 4]{CE}.
\end{proof}

\section{$G$-sheaves} \label{sec-Gsh}
Grothendieck in \cite{Grot} introduces the concept of a $G$-sheaf, where $G$ is a
group.  Note that Grothendieck uses ``sheaf'' to mean a presheaf which satisfies (S1)
and (S2), and uses ``espace \'etal\'e'' for what we call a sheaf.  The following 
coincides with his definition when this difference in terminology is taken into account.

\begin{definition}\rm Let $X$ be a space on which $G$ acts on the left.  A \textit{$G$-sheaf
over
$X$} is a sheaf of $R$-modules $ \pi\co \s{A}\to  X$ over $X$ on which $G$ acts such that, 
if $a \in \s{A}_x$, then $ga \in \s{A}_{gx}$ for all $x \in X$ and $g \in G$. Thus $G$ acts
on $\s{A}$ by fibre homeomorphisms.
\end{definition}

\begin{definition}\rm Let $s\co X \to \s{A}$ be a section of the $G$-sheaf $\s{A}$ and let $g
\in G$. We define the section $gs$ by 
$$ gs(x) = g\cdot s(g^{-1}x)$$
for all $x \in X$.
\end{definition}

Observe that, for $g,g' \in G$, we have
$$(g (g's))(x) = g\cdot (g's)(g^{-1}x) = g \cdot g' \cdot s( {g'}^{-1} g^{-1}x) = (gg')
\cdot s((gg')^{-1}x).$$
We thus have an action of $G$ on the $R$-module $\Gamma_X \s{A}$ of global sections
of $\s{A}$. Note that in the above, no topology on $G$ is assumed.  Requiring that the action be
given by a continuous function $G \times \s{A} \to \s{A}$ would force the action on $\Gamma_X
\s{A}$ of a connected group  to be trivial.

We denote by $I^G$ the functor from $RG$-modules to $R$-modules defined by $I^G(A)
= A^G =\{a \in A \mid ga = a \text{ for all $g \in G$}\}$, the set of $G$-invariants.  We 
denote by $\Gamma^G_X \s{A}$ or $\Gamma^G \s{A}$ the module of $G$-invariant
sections of the $G$-sheaf $\s{A}$.   Thus $\Gamma^G \s{A} = I^G(\Gamma \s{A})$.
 
We now investigate the category of $G$-sheaves over $X$.  A $G$-homomorphism
$f\co \s{A} \to \s{B}$ of $G$-sheaves over $X$ is a homomorphism of sheaves such that for
all $a \in \s{A}$ and $g \in G$, we have $g \cdot f(a) = f(g \cdot a)$.   We clearly have
an $R$-module $\Hom_G(\s{A}, \s{B})$ of $G$-homomorphisms.  

 Note that if $f\co X \to Y$ is a $G$-map of $G$-spaces, then the direct image $f_*\s{A}$ of
a $G$-sheaf $\s{A}$ over $X$ is a $G$-sheaf over $Y$, and the inverse image 
$f^{-1}\s{B}$ of a $G$-sheaf over $Y$ is a $G$-sheaf over $X$.  

\begin{lemma}  Let $\s{A}$ be a $G$-injective sheaf over $X$ and let $f\co X \to Y$ be a
$G$-map.  Then the direct image $f_*\s{A}$ is a $G$-injective sheaf.
\end{lemma}
\begin{proof} The argument of  \cite[Lemma XIV.3.2, p.\ 157]{Bar} applies.
\end{proof}

To study derived functors, we need to know that every $G$-sheaf is embeddable in a
$G$-injective sheaf.  The following is roughly the argument given in Grothendieck and is
a straightforward modification of the standard argument given for ordinary sheaves in
Bredon \cite{Bred2} and  in Barnes \cite{Bar}.

Let $K$ be a not necessarily commutative ring.  First, I consider a natural embedding $i\co A
\to J(A)$ of a $K$-module $A$ in an injective $K$-module $J(A)$.  Such exist.  For
example, put $A^* = \Hom_\Z(A, \Q/\Z)$.  Then $A^*$ is a right $K$-module with action
$fk(a) = f(ka)$.  Note that for every $a \in A$, there exists some $f \in A^*$ for which
$f(a) \ne 0$.  Let $K\langle A^* \rangle$ be the free right $K$-module on the \textit{set}
$A^*$.  We then have a natural $K$-homomorphism $\eta\co K\langle A^* \rangle \to
A^*$.  We can then put $J(A) = ( K\langle A^* \rangle)^*$ and take $i$ the natural
inclusion $A \to A^{**}$ followed by the dual $\eta^*$ of $\eta$.   I now extend this
functor $J$ to $G$-sheaves.

Let $\s{A}$ be a $G$-sheaf over $X$.  The stalk $\s{A}_x$ over $x$ is an
$RG_x$-module where $G_x$ is the stabiliser of $x$.  We embed $\s{A}_x$ in the
injective $RG_x$-module $J(\s{A}_x)$ as above.  Setting $J(\s{A})(U) =  \Pi_{x \in U}
J(\s{A}_x)$ defines a presheaf.  Now let $\s{J}(\s{A})$ be the sheaf of germs of this
presheaf $J(\s{A})$.  This gives a natural embedding  $i\co \s{A} \to \s{J}(\s{A})$ of
$G$-sheaves over $X$ into injective sheaves,   (injectives of the category of sheaves) as
this is a version of the standard construction if we forget the action of $G_x$.

The naturality of $i$ ensures that $\s{J}(\s{A})$ is also a $G$-sheaf.   We want to show
that it is an injective of the category of $G$-sheaves over $X$.  Let $X_d$ be $X$ taken
with the discrete topology. 

\begin{lemma}The sheaf $\s{L} = \Pi_{x \in X} J(\s{A}_x)$ over $X_d$ is $G$-injective.
\end{lemma}
\begin{proof} Given a diagram
\begin{center}
\setlength{\unitlength}{1em}
\begin{picture}(12,4)(-6,-3) 
\put(-6,0){$0$}
\put(0,0){$\s{A}$}
\put(6,0){$\s{B}$}
\put(0,-3.1){$\s{L}$}
\put(-5,0.3){\vector(1,0){4.8}}
\put(1,0.3){\vector(1,0){4.8}}
\put(0.4,-0.4){\vector(0,-1){1.6}}
\multiput(5.7,0)(-1.6,-0.8){2}{\line(-2,-1){1.1}}
\put(2.5,-1.6){\vector(-2,-1){1.3}}
\put(0.7,-1.3){$\phi$}
\put(4.3,-1.3){$\psi$}
\end{picture}
\end{center}
of $G$-sheaves, to construct the required $\psi\co \s{B} \to \s{L}$, we select one $x$ from
each orbit in $X$.  We can then take an $RG_x$-module extension $\psi_x\co \s{B}_x \to
\s{L}_x = J(\s{A}_x)$ of $\phi_x$.  For another point $y = gx$ of the orbit, we define
$\psi_{y}\co \s{B}_{y} \to \s{L}_{y} = J(\s{A}_{y})$ by setting $\psi_{y}(b) = g
\psi_x(g^{-1} b)$.  This is independent of the choice of $g$, for if $y = g'm$, then $g' =
gh$ where $h \in G_x$.  We then have
$$g' \psi_x({g'}^{-1} b) = gh \psi_x({g'}^{-1} b) = g \psi_x(h {g'}^{-1} b) = g
\psi_x(g^{-1} b)$$
as required.
\end{proof}

\begin{lemma} $\s{J}(\s{A})$ is $G$-injective. \end{lemma}
\begin{proof}  As in Barnes \cite[Lemma XIV.3.3, p.\ 157]{Bar}, we have $f_*\s{L} =
\s{J}(\s{A})$ where $f\co X_d \to X$ is the identity on the underlying set.  It follows that
$\s{J}(\s{A})$ is
$G$-injective.
\end{proof}

The sheaf $\s{J}(\s{A})$ is clearly injective.  Thus any $G$-sheaf $\s{A}$ over $X$ can be
embedded in a sheaf  which is both injective and $G$-injective.

\begin{Cor} Let $\s{A}$ be a $G$-injective sheaf.  Then $\s{A}$ is injective.  \end{Cor}

\begin{proof}  $\s{A}$ is a direct summand of $\s{J}(\s{A})$. \end{proof}

It is often convenient to calculate sheaf cohomology using resolutions by flabby sheaves
instead of injective sheaves.  It is clear that the sheaf $\s{F}\s{A}$ of germs of the
presheaf $F(U) = \Pi_{x \in U} \s{A}_x$ of not necessarily continuous sections is flabby
and is a $G$-sheaf.

For any $RG$-module $A$, we have the constant sheaf $\s{C}(A) = A \times X$ over
the $G$-space $X$.  This is a $G$-sheaf with the diagonal action $g(a, x) = (ga, gx)$ for
all $g \in G$, $a \in A$ and $x \in X$.  If $X$ is connected, $\Gamma_X \s{C}(A) = A$.  We
show that $\s{C}$ is a left adjoint to $\Gamma_X$.

\begin{lemma}\label{lem-const} Let $A$ be any $RG$-module and let $A \times X$ be 
the constant sheaf over $X$.  Let $\s{S}$ be a $G$-sheaf over $X$. Then $\Hom_G(A
\times X, \s{S}) \simeq \Hom_{RG}(A, \Gamma_X\s{S})$. This isomorphism is natural 
in both arguments.
\end{lemma}
\begin{proof} The module $A$ can be considered as a sheaf over the one-point space $Y$. 
Consider the map $f\co X \to Y$.  The pullback $f^*A $ is the constant sheaf $A \times X$,
while $f_* \s{S} = \Gamma_X \s{S}$.  Thus the result is a special case of the adjunction
$$\Hom_G(f^*\s{T}, \s{S}) = \Hom_G(\s{T}, f_*\s{S})$$
for a $G$-map $f\co X \to Y$ and $G$-sheaves $\s{S}, \s{T}$ over $X, Y$. 
(It is clear that the adjunction $\Hom(f^*\s{T}, \s{S}) = \Hom(\s{T}, f_*\s{S})$ 
respects the $G$-actions.)
\end{proof}   

From this, we immediately obtain (\cite[Corollaire p.\ 198]{Grot})

\begin{Cor} \label{cor-injsec}  If $\s{A}$ is  $G$-injective, then $\Gamma\s{A}$ is
an injective $RG$-module. \end{Cor}

\begin{proof} $\Hom_G(\_{} \times X, \s{A})$is an exact functor.  By Lemma
\ref{lem-const}, $\Hom_{RG}(\_{}, \Gamma_X\s{A})$ is exact. \end{proof}

We have $\Gamma^G = I^G \circ \Gamma_X$.  By Corollary \ref{cor-injsec}, the
conditions for the composite functor spectral sequence are satisfied.  Thus we have the
theorem of Grothendieck \cite[Th\'eor\`eme 5.2.1]{Grot}:

\begin{theorem}\label{th-ssG} There exists a spectral sequence with target
$R^\bullet\Gamma^G (\s{A})$ and 
$$E_2^{pq} = H^p(G; H^q(X; \s{A})).$$
\end{theorem}

This spectral sequence may be constructed by starting with any
$\Gamma_X$-acyclic resolution $\s{Q}^\bullet$ of $\s{A}$, applying $\Gamma_X$ and
then using any $I^G$-acyclic resolutions of the short exact sequences
$$0 \to Z^q \to \Gamma_X \s{Q}^q \to B^{q-1} \to 0$$
and 
$$ 0 \to B^q \to Z^q \to H^q \to 0.$$
By Barnes \cite[Theorem VII.1.2, p.\ 70]{Bar}, the spectral sequence is independent of
the choices of resolutions from page 2 onwards.  (The  assumption in
\cite{Bar} that the $\Gamma_X \s{Q}^q$ are $I^G$-acyclic is not used anywhere in the
proof.)  The choices of resolutions can be made functorially using exact functors.  This
gives a spectral sequence constructor for $(\Gamma_X, I^G)$.  By \cite[Theorem X.5.4,
p.\ 109]{Bar}, all spectral sequence constructors for  $(\Gamma_X, I^G)$ construct the
same spectral sequence from page 2 onwards.

\section{Equivariant sheaf cohomology} \label{sec-EqSh}

I copy the definition of equivariant (singular) cohomology to define equivariant sheaf
cohomology.

\begin{definition}\rm The \textit{equivariant cohomology} of the $G$-sheaf $\s{A}$ over the
$G$-space $X$ is the sheaf cohomology of the sheaf $\s{A}_G = \s{A} \times_G EG$ over 
the space $X_G = X \times_G EG$, that is, 
$$H^\bullet_G(X; \s{A}) = H^\bullet(X_G; \s{A}_G).$$
\end{definition}

Now consider the group of $G$-invariant sections $\Gamma_X^G \s{A}$.

\begin{lemma} Let $\s{A}$ be an $G$-sheaf over the $G$-space $X$.  Then
$$\Gamma_X^G \s{A} \simeq \Gamma_{X_G} \s{A}_G.$$ \end{lemma}

\begin{proof}  Let $f\co X \to \s{A}$ be a $G$-invariant section of $\s{A}$.  We construct a
section $f_G\co X_G \to \s{A}_G$ by setting $f_G(x \times_G e) = f(x) \times_G e$.  As
$$f(gx) \times_G e = g \cdot f(x) \times_G e = f(x) \times_G g^{-1} e,$$
$f_G$ is well-defined.

Conversely, let $\phi\co X_G \to \s{A}_G$ be a section of $\s{A}_G$.   Choose a point $b
\in BG$ and a point $e \in EG_b$.  Then $x \mapsto x \times_G e$ is a homeomorphism
of $X$ onto the fibre $(X_G)_b$.  Similarly, we have a homeomorphism of $\s{A}$
onto the fibre $(\s{A}_G)_b$ using the same choice of $e$.  The restriction $\phi \vert
(X_G)_b \to (\s{A}_G)_b$ thus defines a section $f\co X \to \s{A}$.  A loop $\gamma$ in
$BG$ at $b$, representing an element  $g$ in the fundamental group $G$ of $BG$ gives a
map of $(X_G)_b$ to itself.  This map is the action of $g$ on $X$.  Similarly, $\gamma$
gives the action of $g$ on $\s{A}$ and it follows that $f(gx) = g \cdot f(x)$ for all $x
\in X$.  Thus $f$ is a $G$-invariant section.

It is clear that these constructions are mutually inverse isomorphisms.
\end{proof}

\begin{lemma}  There exists a natural transformation $\eta\co R^\bullet \Gamma_X^G \to 
H^\bullet_G(X; \_)$.
\end{lemma}

\begin{proof}  Let $0 \to \s{A} \to \s{Q}^0 \to \s{Q}^1 \to \ldots $ be an injective
$G$-resolution of $\s{A}$.  Then $0 \to \s{A}_G \to \s{Q}^0_G \to \s{Q}^1_G \to 
\ldots$ is an exact sequence of sheaves over $X_G$.  Let $0 \to \s{A}_G \to \s{X}^0 \to
\s{X}^1 \to \ldots$ be an injective resolution of $\s{A}_G$ in the category of sheaves
over $X_G$.  Then there exists a translation $t\co \s{Q}^\bullet_G \to \s{X}^\bullet$ over the
identity, and all such are homotopic.   Applying the functor $\Gamma_{X_G}$ and taking
cohomology gives the result.
\end{proof}

If the $\s{Q}^n$ can be chosen such that the $\s{Q}^n_G$ are $\Gamma_{X_G}$-acyclic, then
the resulting homomorphism $R^\bullet \Gamma_X^G (\s{A}) \to H^\bullet_G(X; \s{A})$
is an isomorphism.   In Section \ref{sec-SwL}, I prove the Theorems \ref{th-eacyc},
\ref{th-equivar} that this holds if $G$ is discrete and
$X$ is compact, but it may not hold in general.  Note that except for trivial cases,
$\s{Q}^n_G$ cannot be flabby and so is not injective.  To see this, suppose the stalk
$\s{Q}_x \neq 0$.  We then have sections $s_1 \neq s_2$ over some neighbourhood $U$
of $x$ for which the corresponding $\pi | s_i(U) \to U$ are homeomorphisms.  Now $BG$
is locally contractible, and we can take a point $b \in BG$ and neighbourhood $W$ of
$b$ such that its inverse image in $EG$ is $W \times X$ and such that $W$ is
path-connected.  Then $\s{Q}_G \mid W \times X = W \times \s{Q}$.  We can take
disjoint open sets $W_1, W_2 \subset W$ and a section $s\co (W_1 \cup W_2)
\times U \to \s{Q}$ which, on $W_i \times U$ is $W_i \times s_i$.  This section cannot
be extended to a section over $W \times U$. 

\section{The Swan spectral sequence} \label{sec-Swan}

Sikora defined the spectral sequence he calls the Swan spectral sequence only for a
constant sheaf.  The definition, reproduced in Section \ref{sec-Intro}, is easily
generalised to arbitrary $G$-sheaves.

\begin{definition}\rm  Let $\s{Q}^\bullet$ be a $G$-injective resolution of the $G$-sheaf
$\s{A}$ over $X$.  Let $P_\bullet$ be a projective resolution of $R$ as $RG$-module. 
The \textit{Swan spectral sequence} of $\s{A}$ is the spectral sequence of the double complex
$$\Sw^{pq}(\s{A}) = \Hom_{RG}(P_p, \Gamma_X\s{Q}^q)$$
using $p$ as filtration degree. \end{definition}

\begin{theorem}  Let $P_\bullet$ be a projective resolution of $R$ as $RG$-module.  Let
$\s{Q}^\bullet$ be a $G$-resolution of the $G$-sheaf $\s{A}$ by  $\Gamma_X$-acyclic
sheaves.  Then the spectral sequence, given by the double complex 
$$ \Sw^{pq}(\s{A}) = \Hom_{RG}(P_p, \Gamma_X\s{Q}^q)$$
using $p$ as filtration degree converges with target the right derived functors of
$\Gamma_X^G$.  The spectral sequence has
$$E_1^{pq}(\Sw(\s{A})) = \Hom_{RG}(P_p, H^q(X; \s{A})),$$
$$ E_2^{pq}(\Sw(\s{A})) = H^p(G; H^q(X; \s{A})),$$
and is independent of the choice of resolutions from page $2$ onward.
\end{theorem}

\begin{proof}  The sequence is first quadrant, so converges strongly to the cohomology
of the total complex.  We consider first the case where the $\s{Q}^q$ are
$G$-injective.  We calculate the cohomology of the total complex using the other spectral
sequence
$${E'}_0^{pq} = \Hom_{RG}(P_q, \Gamma_X\s{Q}^p)$$
of the double complex.  By Corollary \ref{cor-injsec}, $\Hom_{RG}(\_, \Gamma_X
\s{Q}^p)$ is an exact functor. Therefore
$${E'}_1^{pq} = \begin{cases} \Gamma^G_X (\s{Q}^p) &\text{for $q= 0$}\\
                                     0 &\text{for $q > 0$.} \end{cases}$$
It follows that 
$${E'}_\infty^{pq} = {E'}_2^{pq} = \begin{cases} R^p\Gamma_X^G(\s{A}) &\text{for $q=
0$}\\
                                     0 &\text{for $q > 0$.} \end{cases}$$
Thus the cohomology of the total complex is $R^\bullet \Gamma_X^G(\s{A})$ as
asserted.

We now calculate $E_1$ and $E_2$.   Since the $P^p $ are projective $RG$-modules, 
$$ E_1^{pq} = \Hom_{RG}(P_p, H^q(\Gamma_X \s{Q}^\bullet)) = \Hom_{RG}(P_p,
H^q(X; \s{A}))$$ 
and
$$ E_2^{pq} = H^p(G; H^q(X; \s{A})).$$

Now consider a $G$-resolution $0 \to \s{A} \to \s{F}^0 \to \s{F}^1 \to \ldots $ by
$\Gamma_X$-acyclic sheaves, and form the double complex 
$$ F_0^{pq} = \Hom_{RG}(P_p, \Gamma_X \s{F}^q).$$
There exists a $G$-translation $t\co \s{F}^\bullet \to \s{Q}^\bullet $.  This induces a
translation $t\co F_0^{\bullet\bullet } \to E_0^{\bullet\bullet}$ of double complexes. 
As the translation of resolutions induces an isomorphism  $H^\bullet (\Gamma_X
\s{F}^\bullet) \to H^\bullet(\Gamma_X \s{Q}^\bullet) = H^\bullet(X; \s{A})$, these
spectral sequences are isomorphic from page $1$.  If  $P'$ is another resolution of $R$, the
resolutions
$P$ and $P'$ are homotopy equivalent giving a homotopy equivalence of the filtered complexes,
with homotopies disturbing the filtration degree by at most $1$.   The spectral sequences are
therefore isomorphic from page $2$ onward.
\end{proof}

In particular, we may use the canonical
flabby resolution.  As this is an \mbox{exact} functor from $G$-sheaves to flabby
$G$-resolutions, this gives a spectral sequence \mbox{constructor} for $(\Gamma_X, I^G)$. 
That it satisfies the total complex and $E^{\bullet 0}_1$ axioms for a constructor follows
from the theorem.  Thus the Swan spectral sequence is the Grothendieck composite functor
spectral sequence considered earlier.  Note also  that the resolutions considered in Theorem
\ref{th-canon} for a $G$-sheaf are $G$-resolutions.  Thus if $X$ is paracompact and
HLC, the Swan spectral sequence may be constructed using any of these resolutions.

Sikora also considers for the constant sheaf $R \times X$, using the singular cochain
complex in place of an injective resolution of $R \times X$.   There is a chain map
from the singular cochain complex to the complex $\Gamma_X \s{S}^\bullet(X;
R \times X)$ which induces the isomorphism in cohomology.  Applying $\Hom(P_\bullet,
\ )$ to this gives a translation of the double complexes.  Thus for a paracompact HLC
space and constant sheaf, use of the singular cochain complex again gives the same
spectral sequence.

\section{The Leray spectral sequence} \label{sec-Leray}
Let $\pi\co E \to B$ be a fibre bundle with fibre $F$ and let $\s{A}$ be a sheaf over $E$.  
We have $\Gamma_E \s{A} = \Gamma_B \pi_* \s{A}$.  The Leray spectral sequence is
the Grothendieck spectral sequence of the composite functor $\Gamma_E = \Gamma_B
\circ \pi_*$.  The usual theory of the composite functor spectral sequence applies, and
we can construct the spectral sequence starting with any $\pi_*$-acyclic resolution
$$0 \to \s{A} \to \s{Q}^0 \to \s{Q}^1 \to \ldots $$
of $\s{A}$.  I want to use the cochain sheaf resolutions of Section \ref{sec-coch}. 
To this end, I consider the direct images under $\pi$ of the cochain sheaves
$\sS^n(E; \s{A})$.  To simplify notation and to make the results generally applicable, I
suppose that we have a sheaf $\s{Q}$ over $E$.

The sheaf $\s{T} = \pi_*\s{Q}$ is defined in terms of the presheaf $T$ over $B$ which
assigns to the open subset $U$ of $B$ the module $T(U) = \s{Q}(\pi^{-1}(U))$.   The
sheaf $\s{T} = \pi_*\s{Q}$ is the sheaf of germs of $T$.   Since $E$ is locally a product, to
determine the stalk over $b \in B$, we need only consider open  neighbourhoods $U$ of
$b$ such that $\pi^{-1}(U) \simeq U \times F$.  So we suppose that $E = B \times F$.

Now consider an element $a$ of the stalk $\s{T}_b$.  Then $a$ is the germ of some $t \in
T(U)$ for some $U$.  But $t \in \s{Q}(\pi^{-1}U)$ is a section of $\s{Q}$ over
$\pi^{-1}U$ and so defines a section $t' = t | E_b$ of the restriction of $\s{Q}$ to $E_b$. 

Conversely, let $t'$ be a section of $\s{Q} | E_b$.  For each point $x \in E_b$, there is a
neighbourhood $N_x$ of $t'(x)$ in $\s{Q}$  homeomorphic to its projection
$W_x  \subseteq E$.   By setting $t_x(e) = N_x \cap \s{Q}_e$ for all $e \in W_x$, we
obtain a section $t_x$ of $\s{Q}$ over $W_x$ with $t_x(x) = t'(x)$.  It follows that $t_x |
(E_b \cap W_x)$ coincides with $t'$ on some neighbourhood of $x$ in $E_b$.  By
replacing $W_x$ by a smaller neighbourhood, we obtain a neighbourhood $W_x = U_x
\times V_x$ and section $t_x$ over $W_x$ with $t_x | V_x = t' | V_x$.  We want to put
these together to get a section over some $U_b \times E_b$.  

\begin{lemma}\label{lem-extend} Let $t'$ be a section of $\s{Q} | E_b$.  Then there exists a
section $t \in \s{Q}(W)$ for some open set $W \supseteq E_b$ such that $t|E_b = t'$.
\end{lemma}

\begin{proof}  For each $x \in E_b$, we have a neighbourhood $W_x = U_x \times V_x$
and a section $t_x$ over $W_x$ as above.  Since $E_b$ is paracompact, there exists a
locally finite refinement $\{V_j \mid j \in J\}$ of the covering of $E_b$ and sections $t_j$
over $U_j \times V_j$ which agree with $t'$ on $V_j$.   This covering $\{V_i\}$ is
shrinkable, that is, there exists an open covering $\{V'_j\}$ with the closure $\bar{V'_j}
\subseteq V_j$.  (See Dugundji \cite[Theorem 6.1, p.152]{Dug}.)  Put $W'_j = U_j
\times \bar{V}_j$.

Replacing $V_x$ by a smaller neighbourhood, we may suppose that $V_x$ meets only
finitely many of the $V_j$.  Suppose $x \in V_j$.  Then $t_x$ agrees with $t_j$ on some
neighbourhood of $x$.  There are only finitely many such $j$, so by taking the intersection
of these neighbourhoods, we may suppose that $t_x$ agrees with $t_j$ on $W_x$ for every
$j$ with $x \in V_j$.  Now suppose that $V_x \cap V_j \not=\emptyset$ but $x \not\in
V_j$.  Replacing $V_x$ by its intersection with the complement of $\bar{V_j}$ ensures that
$V_x \cap  \bar{V_j} = \emptyset$.  We may do this for every such $j$ as there are only
finitely of them.  We now have that for every $j$, $t_x$ agrees with $t_j$ on $W_x \cap
W'_j$.  

Now suppose $e \in W_x \cap W_y$.  Then $e \in W'_j$ for some $j$ and
$t_x(e) = t_j(e) = t_y(e)$.  Since these sections $t_x$ agree on intersections, they assemble
into the required section $t$ over $W = \bigcup_x W_x$. \end{proof}

Now \looseness=-1 suppose that $F$ is compact.  Then $\pi\co E \to B$ is a closed map and every
neighbourhood $W$ of $E_b$ contains a neighbourhood of the form $U \times E_b$.  Thus
every element $a$ of the stalk $(p_*\s{Q})_b$ gives a section of $\s{Q} | E_b$, and 
every section of  $\s{Q} | E_b$ arises in this way.  If $a, a'$ give the same section, then we
have sections $t, t'$ over some $U \times F$ and $U' \times F$ which agree on $E_b$.  They
therefore agree on some neighbourhood of $E_b$.  Using the compactness of $E_b$, we
obtain some $U'' \times E_b$ on which they agree.  Thus they give the same germs, that is,
$a =a'$.

 We have thus proved (\cite[4.17.1]{God}):

\begin{lemma} \label{lem-fcomp} Suppose the fibre $F$ is compact.  Then $\s{T}_b =
\Gamma(\s{Q} | E_b)$. \end{lemma}

To go further, we restrict attention to sheaves which have a bundle structure matching that
of $\pi\co E \to B$.   The bundle can be described in terms of its structure group and
coordinate transformations.  Thus we have the fibre $F$, an effective topological
transformation group $G$ of $F$, a family of coordinate neighbourhoods $\{U_j \mid j
\in J\}$, for each $j \in J$, a homeomorphism $\varphi_j \co U_j \times F \to \pi^{-1}(U_j)$
and, for each pair $i,j \in J$, a continuous map $g_{ji} \co U_i \cap U_j \to G$ giving the
coordinate transformations.  (See Steenrod \cite[\S 2.3, p.\ 7]{Steen}.)  If we have a
$G$-sheaf $\s{F}$ over $F$, we can use these same transition functions to define a sheaf
over $E$.  To allow for cases which arise in equivariant cohomology in which the
given group $G$  acting on $\s{F}$ may not act effectively on $F$, we define a bundle
sheaf a little more generally.  We shall suppose that we have a $G$-sheaf $\s{F}$ over a
$G$-space $F$ and that $G$ acts effectively on $\s{F}$.  Let $K = \{g \in G \mid gx=x
\text{ for all } x\in F\}$, and let $\eta$ be the natural homomorphism $G \to \bar{G} =
G/K$.

\begin{definition}\rm We call a  sheaf $\s{A}$ over the total space $E$ of the fibre bundle
$\pi\co E \to B$ a \textit{bundle sheaf} over $E$ if $\s{A}$ is a $G$-bundle over $B$ with
coordinate transformations $g_{ji}$ such that $\eta  g_{ji}$ are the coordinate
transformations of the
$\bar{G}$-bundle $E$ over $B$. \end{definition}  

Equivalently, $\s{A}$ is called a bundle sheaf over $\pi\co E \to B$ if, for some principal
$G$-bundle $P$ over $B$ and $G$-sheaf over $F$, we have $E = P \times_G F$ and $\s{A}
= P \times_G \s{F}$.  In this situation, we shall refer to $\s{F}$ as the fibre sheaf of $\s{A}$.

We can give an alternative description not directly mentioning the group $G$.  Let $p\co
\s{A} \to E$ be the sheaf projection and let $s_0\co E \to \s{A}$ be the zero section.  Then
$\s{A}$ is a bundle sheaf over $E$ if $\pi p\co \s{A} \to B$ is a fibre bundle and $p, s_0$ are
bundle maps.

Suppose  we are given a $G$-sheaf $\s{F}'$ over $F$ and a $\bar{G}$-bundle $\pi\co E
\to B$ with fibre $F$.  Constructing a bundle sheaf $\s{A}'$ from this involves lifting the
transition functions $\bar{g}_{ji} \co U_i \cap U_j \to \bar{G}$ to functions $g_{ji} \co U_i
\cap U_j \to G$ satisfying all the consistency conditions.  This might not be possible. 
However, if we already have some bundle sheaf $\s{A}$ over $E$, then we have such a lifting. 
We can take the associated principal bundle $P$ of $\s{A}$ and construct the desired bundle
sheaf by setting $\s{A}' =  P \times_G \s{F}'$, that is, by  using the coordinate
transformations of
$\s{A}$.  We denote this sheaf by $\s{F}' \times_{\s{A}} B$.  It is a product of $\s{F}'$ and
$B$ twisted by the transition functions of $\s{A}$.  Note that $\s{F} \times_{\s{A}} B =
\s{A}$ and that $\_ \times_{\s{A}} B = P \times_G \_$ is an exact functor from $G$-sheaves
on $F$ to bundle sheaves over \mbox{$\pi\co E \to B$.}

For a bundle sheaf $\s{A}$, we can obtain the result of Lemma \ref{lem-fcomp} without the
need for the assumption that $F$ is compact.

\begin{lemma} Suppose $\s{A}$ is a bundle sheaf over $\pi\co E \to B$ and that $B$ is locally
path connected.  Then $\pi_*(\s{A})_b = \Gamma_{E_b}(\s{A}|E_b)$. \end{lemma}

\begin{proof}  We can take a path connected neighbourhood $U$ of $b$ contained in some
coordinate neighbourhood.  Then $\s{A} | \pi^{-1}(U) = \s{F} \times U$.  For any path
$\sigma\co I \to U$, the restriction of $\s{A}$ to $(\sigma(I), x) \subset U \times F$ is
constant.  By Lemma \ref{constsec}, any section of it is constant.  It follows that any section of
$\s{A} | E_b$ has a unique extension to a section of $\pi^{-1}(U)$.  \end{proof} 

Investigation of the Leray spectral sequence for a bundle sheaf $\s{A}$ would be easy if we
could show that its fibre sheaf being $\Gamma_F$-acyclic implies that $\s{A}$ is
$\pi_*$-acyclic.  It seems unlikely that this is true in general.  I have only been able to
prove it (Lemma \ref{lem-pifine} below) under the assumption that the fibre $F$ is compact.

\begin{lemma} \label{lem-piinj}  Suppose $\s{Q}$ is fine.  Then $\s{Q}| E_b$ is a fine sheaf
over $E_b$. \end{lemma}
\begin{proof}  Let $\s{V} = \{V_i \mid i \in I\}$ be a locally finite open covering of $E_b$. 
Let $U$ be an open neighbourhood of $b$ in some coordinate neighbourhood. Put
$W_i = U \times V_i$.  Then $V_i = E_b \cap W_i$.   Put $W_0 = E - E_b$.  Then $\s{W} =
\{W_i \mid i \in I \cup \{0\}\}$ is a locally finite open covering of $E$.  The fineness
endomorphisms of  $\s{Q}$ for $\s{W}$ restrict to the required endomorphisms of $\s{Q}|
E_b$ for $\s{V}$. \end{proof}

An immediate consequence of Lemmas \ref{lem-fcomp}, \ref{lem-piinj} is that
$R^q\pi_* \s{A}$ is a sheaf over $B$ with stalks $(R^q\pi_* \s{A})_b = H^q(E_b;
\s{A}|E_b)$.  We denote this sheaf by $\s{H}^q(E_b; \s{A}|E_b)$.  It follows that, if the
fibres are compact, then the Leray spectral sequence has
$$E^{pq}_2 = H^p(B; \s{H}^q(E_b; \s{A}|E_b)).$$

\begin{lemma} \label{lem-pifine}  Let $\s{A}$ be a bundle sheaf over $\pi\co E \to B$ with
fibre sheaf $\s{F}$ over the fibre $F$.  Suppose $F$ is compact and that $\s{F}$ is
$\Gamma_F$-acyclic.  Then $\s{A}$ is $\pi_*$-acyclic. \end{lemma}

\begin{proof}  Take an injective resolution
$ 0 \to \s{A} \to \s{Q}^0 \to \s{Q}^1 \to \ldots $
of $\s{A}$.  Applying $\pi_*$ gives the sequence
\begin{equation} 0 \to \pi_*\s{A} \to \pi_*\s{Q}^0 \to \pi_*\s{Q}^1 \to \ldots
\label{e1}\end{equation}
of sheaves over $B$.  We have to show this sequence is exact.  By Lemma
\ref{lem-fcomp}, the sequence of stalks over $b$ is
\begin{equation} 0 \to \Gamma_{E_b}(\s{A}|E_b) \to \Gamma_{E_b}(\s{Q}^0| E_b) \to
\Gamma_{E_b}(\s{Q}^1| E_b) \to \ldots .\label{e2}\end{equation}

By Lemma \ref{lem-piinj}, 
$ 0 \to \s{A}| E_b \to \s{Q}^0| E_b \to \s{Q}^1| E_b \to \ldots $
is a fine resolution of $\s{A}| E_b$.  But $\s{A}| E_b \simeq \s{F}$ is
$\Gamma_{E_b}$-acyclic.  It follows that the sequence of stalks \eqref{e2} is exact and so
also the sequence \eqref{e1} of sheaves.  Therefore
$\s{A}$ is $\pi_*$-acyclic. \end{proof}

\begin{lemma} \label{lem-pilc}  Suppose $\s{A}$ is a bundle sheaf over $\pi\co E \to B$
with the fibre  $F$ compact.  Then $\pi_*(\s{A})$ is a locally constant sheaf over $B$ with
stalks
$\pi_*(\s{A})_b =
\Gamma_{E_b}(\s{A} | E_b)$. 
\end{lemma}

\begin{proof} Any $b \in B$, has a neighbourhood $U$  such that $\s{A}
| \pi^{-1}(U) = \s{F} \times U$.  Thus $\pi_*(\s{A}) | U = \Gamma(\s{F}) \times U$ by
Lemma \ref{lem-fcomp}. 
\end{proof}

It follows that, for a bundle sheaf $\s{A}$, if $E$ is paracompact and HLC and the fibre $F$
is compact, we may apply any of the cochain sheaf resolutions of Section \ref{sec-coch} to
$\s{F}$ extended as twisted product to a bundle sheaf resolution of $\s{A}$ in the first
stage of the construction of the Leray spectral sequence of $\s{A}$.  This may then be
followed by any of those cochain resolutions in the second stage.  In particular, the Leray
spectral sequence of $\s{A}$ is the spectral sequence of either of the double complexes 
$$\Ler^{pq} = \Gamma_B \s{S}^p(B, \Gamma_{E_b} \s{S}^q(E_b)) \text{\quad or \quad}
\ASLer^{pq} = \Gamma_B \bar{\s{S}}^p(B, \Gamma_{E_b} \bar{\s{S}}^q(E_b)).$$ 
For the Alexander-Spanier version $\ASLer$,  the requirement that $E$ be HLC can be
weakened to $B$ locally path connected.

\section{The Leray-Serre spectral sequence} \label{sec-Serre}
I follow the Dress construction \cite{Dress} of the Leray-Serre spectral sequence as explained
in McCleary \cite[pp.\ 208 -- 212]{McC}.  Let $\pi\co E \to B$ be a fibration.  A singular
$pq$-simplex of the fibration is defined to be a commutative diagram
$$\begin{CD}
\Delta^p \times \Delta^q  @>w>> E \\
@VprVV                                        @VV{\pi}V\\
\Delta^p                            @>>u> B,  \end{CD}$$
where $pr$ is the projection onto the first factor.  The homology version of the spectral
sequence is obtained by defining $S_{pq}$ to be the free abelian group on the set of all
$pq$-simplexes.  This gives a double complex, and using $p$ as filtration degree gives
the Leray-Serre spectral sequence.  Dualising, setting $S^{pq} = \Hom(S_{pq}, R)$ gives the
cohomology version.  I want to generalise this to cohomology of an arbitrary sheaf
$\s{A}$ over $E$.  To simplify, I assume that $\pi\co E \to B$ is a fibre bundle and consider
only open sets of $E$ of the form $U \times V$, where $U$ is an open set of $B$ and $V$
is an open set of the fibre.

For any open set $W$ of $E$, I put 
$$D^{pq}(W) = S^{pq}(W; \s{A}(W)) = \Hom(S_{pq}(W), \s{A}(W)).$$
Here, a $pq$-simplex of $W$ is understood to be a $pq$-simplex of $\pi\co E \to B$ as
above, with $w(\Delta^p \times \Delta^q) \subset W$.  Let $\s{D}^{pq}$ be the sheaf of
germs of the presheaf  $D^{pq}$ and let $\Ser^{pq} = \Gamma_E(\s{D}^{pq})$.  This gives
a double complex.  Using $p$ as filtration degree gives my generalisation of the
Dress-Serre spectral sequence.

Another version, analogous to the Bredon generalisation of singular cohomology, is
obtained by setting $\s{D}_{\fg}^{pq}(\s{A}) = \s{D}^{pq}(\s{R}) \otimes \s{A}$ and
$\Ser_{\fg}^{pq} = \Gamma_E(\s{D}_{\fg}^{pq})$.

By McCleary \cite[Theorem 6.28]{McC}, the homology of the total complex $\sum
S_{pq}$ is $H_\bullet(E)$.

\begin{lemma}\label{lem-dex}  Suppose $E$ is HLC.  Then the sequence
$$0 \to \s{D}^{00} \to \sum_{p+q=1} \s{D}^{pq} \to \sum_{p+q=2} \s{D}^{pq} \to 
\ldots $$ is exact. \end{lemma}
\begin{proof}  The argument of Lemma \ref{lem-HLCsim} applies. \end{proof}

\begin{lemma}\label{lem-dbex}  Suppose $E$ is HLC and that $R$ is noetherian.  Then the
sequence
$$0 \to \s{D}_{\fg}^{00} \to \sum_{p+q=1} \s{D}_{\fg}^{pq} \to \sum_{p+q=2}
\s{D}_{\fg}^{pq} \to \ldots $$ is exact. \end{lemma}
\begin{proof}  The argument of Lemma \ref{lem-exfg} applies. \end{proof}

\begin{theorem} \label{th-htot}  Suppose $E$ is paracompact and HLC and that $R$ is
noetherian.  Then the cohomology of each of the total complexes  $\Ser^{\bullet\bullet}$
and
$\Ser_{\fg}^{\bullet\bullet}$ is $H^\bullet(E; \s{A})$.  \end{theorem}

\begin{proof}  By the argument of Lemma \ref{lem-fine}, the sheaves $\s{D}^{pq}$ and
$\s{D}_{\fg}^{pq}$ are fine.  The result follows by Lemmas \ref{lem-dex}, \ref{lem-dbex}.
\end{proof}

We need to show that these constructions do give us the Dress version of the Leray-Serre
spectral sequence in the special case of a constant sheaf $\s{A} = A \times E$.   As we have
seen before, in this case, we need only consider constant sections in the passage to germs. 
Thus $\s{D}^{pq}$ is the sheaf of germs of the presheaf $S^{pq}(W) = \Hom(S_{pq}(W),
A)$,  and $\s{D}_{\fg}^{pq} = \s{D}^{pq}(\s{R}) \otimes A$.   As $E$ is assumed
paracompact,  it follows as proved in Section \ref{sec-coch} for sheaves of cochains, that
$\Gamma_E(\s{D}^{pq}) = S^{pq}(E; A)/ O^{pq}$, where $O^{pq}$ is the submodule
of locally zero cochains.  We thus have an epimorphism of double complexes
$\eta\co S^{\bullet \bullet}(A) \to \Gamma_E\s{D}^{\bullet \bullet} = \Ser^{\bullet
\bullet}$.   Similarly, we have a map of double complexes $S^{\bullet \bullet}(R) \otimes A
\to \Gamma_E(\s{D}_{\fg}^{\bullet \bullet}) = \Ser_{fg}^{\bullet \bullet}$, but, as pointed
out earlier for the map $S^\bullet(X; R) \otimes A \to \Gamma \s{S}_{\fg}^\bullet(X; A)$,
this need not be surjective.

To show that $\eta$ induces an isomorphism of the spectral sequences, we adapt
the theory of small simplexes (Eilenberg and Steenrod \cite[\S VII.8]{ES}) to
bisimplexes.  To avoid confusion, I denote their map $R$ by $\rho$.  It is convenient to
study the bisimplex
$\sigma = (u, w)$ using only the map
$w\co \Delta^p \times\Delta^q \to E$.  A map $w\co \Delta^p \times\Delta^q \to E$
determines a bisimplex $\sigma = (u, w)$ provided that, for each $x \in \Delta^p$, $\pi
w(x, \Delta^q)$ is a single point, a requirement I shall refer to as the bisimplex condition.

A map $w \co \Delta^p \times\Delta^q \to E$ may be considered as a map $\Delta^p \to
E^{\Delta^q}$.  The maps $\Sd, \rho$ of \cite[\S VII.8]{ES} when applied to this give maps
$\Sd' \co S_{pq} \to S_{pq}$ and $\rho' \co S_{pq}  \to S_{p+1,q}$ such that 
$$\partial' \Sd' = \Sd' \partial'\text{ and }1 - \Sd' = \partial' \rho' + \rho' \partial'.$$
  As $\Sd$ and $\rho$ are
natural,  $\Sd'$ and $\rho'$ commute with the maps \mbox{$\partial'' \co E^{\Delta^q} \to 
E^{\Delta^{q-1}}$} induced by the face maps $\partial_i \co \Delta^{q-1} \to \Delta^q$. 
Thus we have 
$\partial'' \Sd' = \Sd' \partial''$ and $\partial'' \rho' + \rho' \partial'' = 0$  (because we
put in $(-1)^p$ in  $\partial''$).  Likewise, regarding $w$ as a map $\Delta^q \to
E^{\Delta^p}$, we obtain $\Sd'' \co S_{pq} \to S_{pq}$ and $\rho''\co S_{pq} \to S_{p,q+1}$,
(giving $\rho''$ the sign $(-1)^p$) such that 
$$\partial'' \Sd'' = \Sd'' \partial''\text{ and }1 - \Sd'' = \partial'' \rho'' + \rho''
\partial''$$
 and also satisfying 
$$\partial' \Sd'' = \Sd'' \partial'\text{ and }\partial' \rho'' + \rho'' \partial' = 0.$$
 
Note that these maps all preserve the bisimplex condition.  Put $\Sd = \Sd' \Sd''$ and
$\rho = \rho' \Sd'' + \rho''$.  For these, we have 
\begin{equation*}
\begin{split}
\partial \rho + \rho \partial &= (\partial' + \partial'')(\rho' \Sd'' + \rho'') +   (\rho' \Sd''
+ \rho'')(\partial' + \partial'')\\
&= \quad \partial' \rho' \Sd'' +\  \partial'' \rho' \Sd'' +\  \partial'\rho''+\partial'' \rho''\\ 
&\quad + \rho' \Sd'' \partial' + \rho' \Sd'' \partial'' + \rho'' \partial' + \rho'' \partial''\\ 
&= (\partial' \rho' + \rho' \partial') \Sd'' + \rho'(-\partial'' \Sd'' + \Sd'' \partial'') + 
(\partial'' \rho'' + \rho'' \partial'') \\
 &= (1- \Sd')\Sd'' + 1 - \Sd'' = 1 - \Sd' \Sd''\\ &= 1 - \Sd
\end{split} \end{equation*}
We thus have the double barycentric subdivision chain map $\Sd \co S_{pq} \to S_{pq}$ and
the chain homotopy $\rho \co S_{pq} \to S_{p+1, q} + S_{p, q+1}$, $1 -\Sd = \partial \rho + 
\rho \partial$.

Now let $\s{W}$ be a covering of $E$ by open sets $W$.  Let $$S_{pq}(E, \s{W}) = \sum
\{S_{pq}(W) \mid W \in \s{W}\}$$ be the subgroup of $S_{pq}(E)$ spanned by the ``small''
bisimplexes,  those whose image is contained in some $W \in \s{W}$.  The
inclusion
$\inc\co S_{pq}(E, \s{W}) \to S_{pq}(E)$ induces the restriction map $\eta \co S^{pq}(E; A) \to
S^{pq}(E, \s{W} ; A)$.

\begin{lemma}\label{lem-bsub} There exists a chain map $\tau\co S_{pq}(E) \to S_{pq}(E,
\s{W})$ and a homotopy $D\co S_{pq}(E) \to S_{p+1,q} + S_{p, q+1}$ such that $\tau
\inc = 1$ and $\partial D + D \partial = 1 - \inc \tau$. \end{lemma}

\begin{proof}  Let $\sigma$ be a $pq$-simplex of $E$.  As in the standard small simplex
theory, we may use the compactness of $\Delta^p\times \Delta^q$ to conclude that $\Sd^n
\sigma \in S_{pq}(E, \s{W})$ for some $n$.  Let $n(\sigma)$ be the least such
$n$.  Then for the faces of $\sigma$, we have $n(\partial'_i \sigma) \le n(\sigma)$ and
$n(\partial''_j \sigma) \le n(\sigma)$.  Put
$$\tau \sigma = \Sd^{n(\sigma)} \sigma + \sum_{i=0}^p (-1)^i \sum_{j =
n(\partial'_i \sigma)}^{n(\sigma) - 1} \rho \Sd^j(\partial'_i \sigma) + (-1)^p \sum_{i=0}^q
(-1)^i
\sum_{j=n(\partial''_i \sigma)}^{n(\sigma) - 1} \rho \Sd^j(\partial''_i \sigma)$$
and put $D \sigma = \sum_{j = 0}^{n(\sigma)-1} \rho  \Sd^j \sigma$.  Then
\begin{equation*}
\begin{split}
\partial D \sigma &=\sum_{j=0}^{n(\sigma)-1} \partial \rho \Sd^j \sigma\\
&= \sum_{j=0}^{n(\sigma)-1}(1 - \Sd - \rho \partial) \Sd^j \sigma  \\
&= \sigma - \Sd^{n(\sigma)} \sigma - \sum_{j=0}^{n(\sigma)-1} \rho \Sd^j \partial \sigma \\
&= \sigma - \Sd^{n(\sigma)} \sigma - \sum_{i=0}^p (-1)^i \sum_{j=0}^{n(\sigma)-1} \rho 
\Sd^j \partial'_i \sigma  -  \sum_{i=0}^q (-1)^{p+i} \sum_{j=0}^{n(\sigma)-1} \rho \Sd^j
\partial''_i \sigma \\
D \partial \sigma &=  \sum_{i=0}^p (-1)^i \sum_{j=0}^{n(\partial'_i\sigma)-1} \rho 
\Sd^j \partial'_i \sigma  +  \sum_{i=0}^q (-1)^{p+i} \sum_{j=0}^{n(\partial''_i\sigma)-1} \rho
\Sd^j \partial''_i \sigma \\
\end{split}
\end{equation*}
Thus $(\partial D + D \partial) \sigma = (1 - \tau) \sigma$.
\end{proof}
\begin{theorem} \label{th-DSconst} Suppose $E$ is paracompact and HLC and  let
$\s{A}$ be the constant sheaf $\s{A} = A \times E$.  Then the spectral sequences of 
$S^{\bullet \bullet}(A)$ and $ \Ser^{\bullet\bullet}(A)$ are isomorphic from page 2 
onwards.
\end{theorem}

\begin{proof}  The $S^{pq}(E, \s{W}; A)$ for open coverings $\s{W}$ of $E$ form a direct
system with restriction maps $\eta_{21}\co S^{pq}(E, \s{W}_1; A) \to S^{pq}(E, \s{W}_2; A)$
whenever $\s{W}_2$ is a refinement of $\s{W}_1$.  The direct limit of this system is
$\Ser^{pq}(A)$.  The maps commute with the differentials and so, with the construction of
the spectral sequences.  Thus the spectral sequence of the direct limit is the direct limit of
the spectral sequences.  By Lemma \ref{lem-bsub} and Cartan and Eilenberg \cite[Prop.
XV.3.1, p.321]{CE},  the $E_2^{pq}(\eta_{ij})$ are isomorphisms, and it follows that the map
$E_2(\eta) \co E_2(S^{\bullet \bullet}(A)) \to E_2(\Ser^{\bullet \bullet}(A))$ is an
isomorphism.  
\end{proof}

For finitely generated constant coefficients $A$, we have $\Ser_{\fg}^{\bullet \bullet}(A) =
\Ser^{\bullet \bullet}(A)$, so the spectral sequences of both $\Ser_{\fg}^{\bullet
\bullet}(\s{A})$ and $\Ser^{\bullet \bullet}(\s{A})$  may be considered to be
generalisations of the Leray-Serre spectral sequence.  The effect on the cohomology of a cochain
complex, of restricting to cochains having values in a finitely generated subgroup
disappears when we work with the complex of sections of the sheaf of germs.  This
suggests that $\Ser_{\fg}^{\bullet \bullet}(\s{A})$ and $\Ser^{\bullet \bullet}(\s{A})$
should give the same spectral sequence from page 2, but I have not been able to prove this.

To get useful results on more general sheaves, I need a structure on the bundle
analogous to the connections used in differential geometry.   I define this in terms of
coordinate neighbourhoods and transition functions.  So suppose $G$ acts effectively on
the fibre $F$ and that the bundle has coordinate functions $\varphi_j \co U_j \times F \to
\pi^{-1}(U_j)$ and transition functions $g_{ji} \co U_i \cap U_j \to G$.  Thus, for $u \in U_i
\cap U_j$,
$$g_{ij}(u) = \varphi^{-1}_j(u, \_) \varphi_i(u, \_) \co F \to F.$$

\begin{definition}\rm A subset $X \subset \pi^{-1}U_j$ is said to be \textit{level} with
respect to
$\varphi_j$ if $X \subseteq \varphi_j(U_j, t)$ for some $t \in F$.  A map $f\co \Delta \to
\pi^{-1}U_j$ is called \textit{level} with respect to $\varphi_j$ if $f(\Delta)$ is level.
\end{definition}

\begin{definition}\rm The family $\{\varphi_j \mid j \in J\}$ is called a \textit{connection}
if, for all
$i, j \in J$ and all path-connected subsets $X$ of $\pi^{-1}(U_i \cap U_j)$, $X$ level with
respect to $\varphi_i$ implies that $X$ is also  level with respect to $\varphi_j$.
\end{definition}

There is a special  case in which the existence of a connection is clear.

\begin{lemma} \label{lem-conn} Suppose that $B$ is locally path connected.  Then the
following conditions are equivalent: \rm
\begin{enumerate}
\item \textit{$E$ has a connection;}
\item \textit{$E$ can be defined using locally constant transition functions;}
\item \textit{$E$ can be defined using a discrete group $G$.}
\end{enumerate} 
\end{lemma}

\begin{proof} $\text{(1)} \Rightarrow \text{(2)}$ We may suppose that $G$ acts
effectively on $F$ and let $\varphi_i, g_{ji}$ define a connection.  Let $X$ be a path
connected subset of $U_i \cap U_j$.  Then for $v \in F$, $Y = \varphi_i(X, v)$ is level in
$\pi^{-1}(U_i)$ and therefore also in $\pi^{-1}(U_j)$.  But $Y = \{\varphi_j(x, g_{ji}(x)v) \mid
x \in X\}$, so
$g_{ji}(x)v$ is independent of $x$ for each $v \in F$.  Since $G$ is effective, $g_{ji}(x)$ is
independent of $x$ for $x \in X$.

$\text{(2)} \Rightarrow \text{(3)}$  The $g_{ji}$, being locally constant, are still continuous
if we give $G$ the discrete topology.

$\text{(3)} \Rightarrow \text{(1)}$
 Let $f\co I \to \pi^{-1}(U_i \cap U_j)$ be level with respect to $\varphi_j$ and  put
\mbox{$u = \pi f \co I \to U_i \cap U_j$}. Then $f(t) = \varphi_j (u(t), c)$ for some $c \in F$
and all $t
\in I$.  So $f(t) = \varphi_i (u(t),  g_{ij}(u(t))(c))$.  But $g_{ij} \circ u \co I \to G$ is continuous
and $G$ is discrete.  Therefore $ g_{ij}(u(t)) = g \in G$ is independent of $t$
and we have $f(t) = \varphi_i (u(t), gc)$.  Thus $f$ is level with respect to $\varphi_i$.  It follows
that any path connected subset of $\pi^{-1}(U_i \cap U_j)$ which is level with respect to 
$\varphi_j$ is also level with respect to  $\varphi_i$.  Thus $\{\varphi_j \mid j \in J \}$ is a
connection.
\end{proof}

\begin{definition}\rm  A $pq$-simplex $(u,w)$ of an open set $W \subseteq \pi^{-1}U_j$ is
called \textit{level} (with respect to $\varphi_j$) if, for all $y \in \Delta^q$, $w(\_, y) \co
\Delta^p \to W$ is level. \end{definition}

If we identify $\pi^{-1}(U_j)$ with $U_j \times F$ via $\varphi_j$, we can simplify the 
notation.  A $pq$-simplex $\sigma = (u, w)$ can be expressed as $w(x,y) = (u(x),
v(x,y))$ for $(x,y) \in \Delta^p \times \Delta^q $, where $v(x,y) \in F$.  The condition
for $\sigma$ to be level becomes that $v$ be a function of $y$ only, thus $w(x,y) =
(u(x), v(y))$.

Let $c \in S^{pq}(W; \s{A}(W))$.  Let $\rho c$ be the restriction of $c$ to level
$pq$-simplexes.  We call these functions $\rho c$ level cochains.  Since all face
operators send level bisimplexes to level bisimplexes, $\rho$ is a map of double
complexes.  Suppose we have a connection, so the concept of level is unambiguous for
small bisimplexes.  Let $\s{N}^{pq}$ be the sheaf of germs of the $\rho c$.  We then have
a map $\rho\co \s{D}^{\bullet \bullet} \to \s{N}^{\bullet \bullet}$.  Setting $N^{pq} =
\Gamma_E(\s{N}^{pq})$, we have a map of double complexes $\rho\co \Ser^{\bullet \bullet}
\to  N^{\bullet \bullet}$.  Restricting to cochains with values in a finitely generated
submodule, we obtain maps $\rho'\co \s{D}_{\fg}^{\bullet \bullet} \to \s{N}_{\fg}^{\bullet
\bullet}$ and $\rho' \co \Ser_{\fg}^{\bullet \bullet} \to N_{\fg}^{\bullet \bullet}$.

\begin{theorem} \label{th-level} Suppose the coordinate functions $\{\varphi_j \mid j \in J\}$
define a connection.   Then the maps $\rho_1\co E^{\bullet \bullet}_1(\Ser) \to E^{\bullet
\bullet}_1(N)$ and $\rho'_1\co E^{\bullet \bullet}_1(\Ser_{\fg}) \to E^{\bullet
\bullet}_1(N_{\fg})$ induced by the restriction $\rho$ to level bisimplexes are
isomorphisms.
\end{theorem}

\begin{proof}  We show that $\rho\co \Ser^{\bullet \bullet} \to  N^{\bullet \bullet}$ and
$\rho'\co \Ser_{\fg}^{\bullet \bullet} \to  N_{\fg}^{\bullet \bullet}$ are chain equivalences
with respect to the second differentials of the double complexes.  Let $\sigma$ be the
\mbox{$pq$-simplex} of $\pi^{-1} U_i$ defined by the map  $(u_i, v_i)\co \Delta^p \times
\Delta^q \to U_i \times F$.  Let $l \sigma$ given by   $(u_i, lv_i)$ be the unique level
$pq$-simplex with $lv_i(0,y) = v_i(0,y)$ for all $y \in \Delta^q$, that is,
$lv_i(x,y) =  v_i(0,y)$.  The map $\gamma\co I \times \Delta^p \times \Delta^q \to
\pi^{-1} U_i$ defined by 
$$\gamma(t,x,y) = \varphi_i(u_i(x), v_i(tx, y))$$
is a homotopy from the level $pq$-simplex $l \sigma$ to $\sigma$.  This homotopy
$\gamma$ does not depend on the neighbourhood $U_i$, for suppose $(u_j, v_j) \co
\Delta^p \times \Delta^q \to U_j \times F$  defines the same map into $E$.  Then $u_i =
u_j$ and 
$$v_j(x,y) = g_{ji}(u_i(x)) v_i(x,y).$$
Defining homotopies by $\gamma_i (t,x,y) = (u_i(x), v_i(tx,y))$ and $\gamma_j$ similarly
for $U_j$, we want to show that these define the same map $I \times \Delta^p \times
\Delta^q \to E$.  For this, we need $v_j(tx, y) = g_{ji}(u_i(x)) v_i(tx, y)$.  This holds as 
$g_{ji}(u_i(tx)) = g_{ji}(u_i(x))$ because the transition functions are locally constant. 

For fixed $p$, the cochains $c$ whose value on $(u, w)$ depends only on $u$ and $w(0,
\_)$ form a subcomplex with respect to the second differential.  This subcomplex can be
identified with $\rho S^{p \bullet}$.  The homotopy $\gamma$ shows it to be a
deformation retract.  This relationship is preserved under the taking of germs and sections. 
Thus the maps $\rho\co \Ser^{p \bullet} \to N^{p \bullet}$ and $\rho'\co \Ser_{\fg}^{p \bullet}
\to N_{\fg}^{p \bullet}$ are chain equivalences and the result follows.
\end{proof}

\section{Comparison of Leray-Serre and Leray} \label{sec-SerL}
To relate the Leray-Serre spectral sequence to the Leray, I construct a map $\psi\co N^{\bullet
\bullet} \to \Ler^{\bullet \bullet}$.  I begin with a lemma which strengthens part of the
standard proof that $S^n(X; \Z) \to \Gamma \s{S}^n(X; \Z)$ is surjective.  Let $\s{A}$ be a
sheaf over the space $X$.  Let $\Delta$ be a test space.  We refer to any map $\sigma \co
\Delta \to U \subseteq X$ as a simplex of $U$ and denote the set of all such by $S(U)$.  A
local cochain is any function $c \co S(U) \to \s{A}(U)$ for some open set $U$.  We say that
the local cochains $c_1, c_2$ defined on open sets $U_1, U_2$ agree on $U_1 \cap U_2$
if, for every simplex $\sigma$ of $U_1 \cap U_2$, $c_1(\sigma)$ and $c_2(\sigma)$ have
the same restriction in $\s{A}(U_1 \cap U_2)$.

\begin{lemma}\label{lem-agree}  Let $\s{A}$ be a sheaf over the paracompact 
space $X$.  Let $\s{S}$ be the sheaf of germs of local cochains and let $\gamma$ be a
section of $\s{S}$.  Then there exists, for each $x \in X$, a neighbourhood $V_x$ of $x$
and a local cochain $c_x$ defined over $V_x$ such that \rm
\begin{enumerate}
\item \textit{the germ of $c_x$ at $y \in V_x$ is $\gamma(y)$, and}
\item \textit{for all $x,y \in X$, $c_x$ and $c_y$ agree on $V_x \cap V_y$.}
\end{enumerate} \end{lemma}

\begin{proof} For $x \in X$, $\gamma(x)$ is the germ at $x$ of some local cochain $c_x$
defined over some neighbourhood $V_x$ of $x$.  Taking the germ at $y \in V_x$ of $c_x$
gives a section of $\s{S}$ over $V_x$ which agrees with $\gamma$ at $x$.  These section
agree on some neighbourhood $V'_x$ of $x$.  Replacing $V_x$ with $V'_x$ gives property
(1).

The covering $\s{V} = \{V_x \mid x \in X\}$ has a locally finite refinement $\s{U} = \{U_j
\mid j \in J\}$.  For each $j$, we can choose a local cochain $C_j$ defined over $U_j$
which is the restriction of some $c_x$.  The covering $\s{U}$ is shrinkable.  Let $\s{U}' =
\{U'_j \mid j \in J\}$ be a shrinking of $\s{U}$.  Consider $x \in X$.  Then $x$ has a
neighbourhood meeting only finitely many of the $U_j$.  By replacing $V_x$ with its
intersection with such a neighbourhood, we may suppose that $V_x$ meets only finitely
many of the $U_j$.

Suppose $x \in U_j$.  Then $c_x$ and $C_j$ have the same germ at $x$, so agree on some
neighbourhood $V_{x,j}$ of $x$.  By replacing $V_x$ with the intersection of these
$V_{x,j}$, we obtain $V_x$ such that $c_x$ is the restriction to $V_x$ of $C_j$ for every
such $j$.  Now suppose $V_x$ meets $U_j$ but $x \notin U_j$.  By replacing $V_x$ by its
intersection with the complement of $\bar{U}'_j$, we may assume that $V_x$ does not
meet $U'_j$.  We may assume this for every such $j$ as there are only finitely many of
them.

We now have $V_x$ such that, if $x \in U_j$, then $V_x \subseteq U_j$ and $c_x =
C_j|V_x$, and, if $x \notin U_j$, then $V_x \cap U_j' = \emptyset$.  Suppose $\sigma$ is
a simplex of $V_x \cap V_y$.  Then there exists $z \in V_x \cap V_y$ and $z \in U'_j$ for
some $j$.  Since $V_x \cap U'_j \ne \emptyset$, $x \in U_j$, $V_x \subseteq U_j$ and
$c_x = C_j|V_x$.  Similarly, $c_y = C_j | V_y$.  Thus $c_x | V_x \cap V_y = c_y |  V_x \cap
V_y$. 
\end{proof}
 
Now suppose that the bundle $\pi\co E \to B$ has $E$ paracompact and HLC, with
compact fibre $F$.  Suppose further, that it has a connection, so we have a family of
coordinate patches $\varphi_i \co U_i \times F \to \pi^{-1}(U_i)$ with locally constant
transition functions.   Let $\s{A}$ be a bundle sheaf over $E$ with fibre sheaf $\s{F}$ over
$F$.  Corresponding to level subsets of $U_i \times F$, we have level sections of $U_i
\times \s{F}$.  An element of a stalk over $(u, v) \in U_i \times F$ is an element of
$\s{F}_v$.  A section over a level path in $U_i \times F$ gives a continuous function from
$I$ into $\s{F}_v$.  As $\s{F}_v$ is discrete, this function must be constant.  It follows that
elements of $\s{A}(\varphi_i (U_i \times F))$ are constant over path connected level subsets. 
As  $B$ is locally path connected, in the passage to germs of level cochains, we need only
consider cochains whose values are level sections.  A level section over $W = \varphi_i(U
\times V)$ is a section of $\s{F}$ over $V$.

An element of $N^{pq}$ is a section $\gamma$ of the sheaf of germs of level
$pq$-cochains.  By Lemma \ref{lem-agree}, we have, for each $e \in E$, a neighbourhood
$W_e$ of $e$ and a level $pq$-cochain $c_e$ defined on $W_e$ with germ
$\gamma(e')$ at $e' \in W_e$ and with the $c_e$ agreeing on overlaps.  We may take
$W_e$ within a coordinate patch, $W_e = \varphi_i(U_e \times V_e)$ and may also assume
$U_e$ path connected.  Further, $c_e$, having values in level sections, has values in
$\s{F}(V_e)$.  Consider the fibre $E_b$ for $b \in U_i$.  The $V_e$ for $e \in E_b$ cover
$F$ which is compact.  We have a finite subcovering $V_{e_1}, \ldots, V_{e_k}$.  Put $U_b
= U_{e_1} \cap \ldots \cap  U_{e_k}$.  For each $p$-simplex $u$ of $U_b$ and
$q$-simplex $v$ of $V_{e_j}$, we have $c_{e_j}(u,v) \in \s{F}(V_{e_j})$.  As the $c_e$
agree on overlaps, we obtain, passing to germs,  for each $u \in U_b$ a section of
$\s{S}^q(E_b; \s{F})$.  Thus we have a $p$-cochain on $U_b$ with values in the locally
constant sheaf $\{\Gamma_{E_b}\s{S}^q(E_b; \s{F})\}$.  Passing to germs gives an element
$\psi(\gamma) \in \Ler^{pq}$.

The element $\psi(\gamma)$ does not depend on the choices made in its construction. 
Suppose for each $e$, we take another level cochain $c'_e$ defined on a neighbourhood
$W'_e$.  Since $c_e, c'_e$ have the same germ at $e$, they agree on some neighbourhood
$W''_e$.  The $W''_e$ for $e \in E_b$ provide a common refinement  of the two coverings
of $E_b$ and, using the compactness, we obtain a finite common refinement.  Passing to
germs, we obtain, for $p$-simplexes $u$ of a sufficiently small neighbourhood of $b$,  the
same section of $\s{S}^q(E_b; \s{F})$ and so the same element of $\Ler^{pq}$.

Denote the first and second differentials of the double complexes by $\partial', \partial''$. 
To calculate $\partial' \gamma$ and $\partial'' \gamma$, we  calculate $\partial' c_e$ and
$\partial'' c_e$ for the $c_e$ used above in the calculation of $\psi(\gamma)$.  We
are working with germs, so we need only consider level $(p+1,q)$- and $(p,
q+1)$-simplexes of the $W_e$.  All faces of these lie in $W_e$ and, as $\psi$ does not
depend on the choices of neighbourhoods, it follows that $\psi$ commutes with $\partial'$
and $\partial''$.  Thus $\psi$ is a well-defined map of double complexes.

Now suppose $\gamma \in N_{\fg}^{pq}$.  Then the cochains $c_e$ give, for each
$p$-simplex $u$ of $U_b$ and $q$-simplex $v$ of $V_{e_i}$, an element $c_{e_i}(u,v)$ of some
finitely generated subgroup $F_i \le \s{F}(V_{e_i})$.  Fixing $u$, we have a $q$-cochain
$c^i_u \in S^q(V_{e_i}; \R) \otimes \s{F}(V_{e_i})$.  Thus $\psi(\gamma) \in \Gamma_B
\s{S}^p(B;
\Gamma \s{S}_{\fg}^q(E_b; \s{F}))$.  Note that, although the values of $c^i_u(v)$ lie in a
finitely generated subgroup of $\s{F}(V_{e_i})$, we have infinitely many $c^i_u$ and these
need not give germs in a finitely generated subgroup of $\Gamma \s{S}_{\fg}^q(E_b;
\s{F}))$.  Thus $\psi(\gamma)$ need not lie in $\Gamma_B \s{S}_{\fg}^p(B; \Gamma
\s{S}_{\fg}^q(E_b; \s{F}))$.

Put ${\Ler'}^{pq} = \Gamma_B \s{S}^p(B; \Gamma \s{S}_{\fg}^q(E_b; \s{F}))$ and
${\Ler''}^{pq} = \Gamma_B \s{S}_{\fg}^p(B; \Gamma \s{S}_{\fg}^q(E_b;
\s{F}))$.  We have the map of double complexes $\psi' \co N_{\fg} \to \Ler'$.  I want
information on the induced map
${\psi'}^{pq}_1 \co E^{pq}_1(N_{\fg}) \to E^{pq}_1(\Ler')$.  I start with $q=0$.

In the construction of $\psi(\gamma)$ for $\gamma \in N^{p0}$, we used neighbourhoods
$U_b, V_{e_i}$ and level cochains $c_e$ defined on $U_b \times V_{e_i}$ with values in
$\s{F}(V_{e_i})$.  We may choose these neighbourhoods path connected.  For each
$p$-simplex $\sigma$ of $U_b$ and each point $v \in V_i$, we get an element of
$\s{F}(V_{e_i})$.  The element $\gamma \in N^{p0}$ is a cocycle with respect to the second
differential $\partial''$ if this element of $\s{F}(V_{e_i})$ is independent of $v$.  Thus
$c_e \in S^p(U_b; \s{F}(V_{e_i}))$.  As the $c_e$ agree on overlaps, we have an element of
$S^p(U_b; \Gamma_F\s{F})$ and $\gamma$ can be identified with the element
$\psi(\gamma) \in E^{p0}_1(\Ler) = \Gamma_B\s{S}^p(B; \Gamma_{E_b}\s{F})$.  Observe
that every element of $\Gamma_B\s{S}^p(B; \Gamma_{E_b}\s{F})$ arises in this way. 
Given a $p$-cochain $c$ on a small neighbourhood of $b \in B$ with values in
$\Gamma\s{F}$, we can make it a level $p0$-cochain by defining $c(\sigma, e) =
c(\sigma)$ for all $e \in E_b$.  Taking germs gives an element of $E^{p0}_1(N)$.  Thus
$\psi^{p0}_1 \co E^{p0}_1(N) \to E^{p0}_1(\Ler)$ is an isomorphism.  

Now suppose that $\gamma \in N_{\fg}^{p0}$. 
Then $c_e$ has values in some finitely generated submodule $A_i$ of $\s{F}(V_{e_i})$. 
The resulting sections of $\s{F}$ are in a submodule of the finitely generated module
$\Pi_iA_i$.  Since $R$ is noetherian, this submodule is finitely generated.   We thus have
a finitely generated module of sections,  a finitely generated  submodule of
$\Gamma_{E_b}\s{F}$.  Thus we have an element of $\s{C}^p(B; \Gamma_{E_b}\s{F})$
and it follows that $\psi(\gamma) \in E^{p0}_1(\Ler'')$.  Going backwards, starting from an
\mbox{element} of $E^{p0}_1(\Ler'')$, we get cochains with values in a finitely generated
submodule of $\Gamma_{E_b}\s{F}$ and so get an element of $N_{\fg}^{p0}$.  Thus
${\psi'}^{p0}_1$ maps $E^{p0}_1(N_{\fg})$ isomorphically onto the submodule
$E^{p0}_1(\Ler'')$ of $E^{p0}_1(\Ler')$.

To get information on ${\psi'}^{pq}_1$ for $q > 0$, we use dimension shifting.  For this,
we need some preliminaries.

\begin{lemma} $N_{\fg}^{p \bullet}, {\Ler'}^{p \bullet}$ and ${\Ler''}^{p\bullet}$ are exact
functors from $G$-sheaves on $F$ to cochain complexes. \end{lemma}

\begin{proof}  $\s{S}_{\fg}^\bullet(F; \_)$ is an exact functor from sheaves on $F$ to
cochain complexes of fine sheaves on $F$.  As fine sheaves on $F$ are
$\Gamma_F$-acyclic, $\{\Gamma_{E_b} \s{S}_{\fg}^\bullet(E_b; \_)\}$ is an exact functor
from $G$-sheaves on $F$ to cochain complexes of locally constant sheaves on $B$.  Thus
$\s{S}^p(B; \Gamma_{E_b} \s{S}_{\fg}^\bullet(E_b; \_)$ and $\s{S}_{\fg}^p(B;
\Gamma_{E_b} \s{S}_{\fg}^\bullet(E_b; \_)$ are exact functors to cochain complexes of 
fine sheaves on $B$.   Taking sections over $B$, it follows that ${\Ler'}^{p \bullet}$ and
${\Ler''}^{p \bullet}$ are exact functors to cochain complexes.

The functor from $G$-sheaves over $F$ to bundle sheaves over $E$ is exact.  The functor
from bundle sheaves to level elements of $\s{S}_{\fg}^{pq}$ is exact.  From a short exact
sequence $0 \to \s{A} \to \s{B} \to \s{C} \to 0$ of bundle sheaves, we get a short exact
sequence $0 \to \s{N}_{\fg}^{pq}(\s{A}) \to \s{N}_{\fg}^{pq}(\s{B}) \to
\s{N}_{\fg}^{pq}(\s{C}) \to 0$.  Since $\s{N}_{\fg}^{pq}(\s{A})$ is fine and therefore
$\Gamma_E$-acyclic, the sequence $0 \to \Gamma_E\s{N}_{\fg}^{pq}(\s{A}) \to 
\Gamma_E\s{N}_{\fg}^{pq}(\s{B}) \to  \Gamma_E\s{N}_{\fg}^{pq}(\s{C}) \to 0$ is also
exact.
\end{proof}

\begin{definition}\rm Let $(\s{C}^\bullet, d)$ be a cochain complex of sheaves over $X$.   We
say that $(\s{C}^\bullet, d)$ is \textit{fine} if, for every locally finite
covering$\{U_\alpha\}$, there exist endomorphisms $\theta_\alpha$ which commute with the
coboundary operator $d$ and satisfy $|\theta_\alpha| \subseteq \bar{U}_\alpha$ and $\sum
\theta_\alpha = 1$.
\end{definition}

\begin{lemma}
It is sufficient to check fineness for a cofinal set of locally finite coverings.
\end{lemma}

\begin{proof} If $\{V_\beta\}$ refines $\{U_\alpha\}$ and we have endomorphisms $\theta_\beta$
for
$\{V_\beta\}$, choose $i(\beta)$ with $V_\beta\subseteq U_{i(\beta)}$
and define $\lambda_\alpha = \sum \theta_\beta$ over the $\beta$ with
$i(\beta)=\alpha$.
\end{proof} 

\begin{lemma}\label{lem-fincomp}
If $X$ is paracompact, the complex $(\s{C}^\bullet,d)$ is fine and $H^n(\s{C}^\bullet)=0$,
then $H^n(\Gamma(\s{C}^\bullet))=0$.
\end{lemma}

\begin{proof}  Suppose $\gamma \in \Gamma \s{C}^n$ and $d \gamma = 0$.  Then, for
each $x \in X$, $d \gamma(x) = 0$.  Since $H^n(\s{C}^\bullet) = 0$, there exists $c_x \in
\s{C}^{n-1}_x$ with $d c_x = \gamma(x)$.  There exists a neighbourhood $N_x$ of $x$
and a section $s_x$ over $N_x$ with $s_x(x) = c_x$.  The section $d s_x$ of $\s{C}^n$
agrees with $\gamma$ at $x$ and so on some neighbourhood $N'_x$ of $x$.  As $X$ is
paracompact, we can take a locally finite refinement $U_\alpha$ of the covering
$\{N'_x\}$ with sections $s_\alpha$ such that $d s_\alpha = \gamma | U_\alpha$.

Let $\{V_\alpha\}$ be a shrinking of $\{U_\alpha\}$ and let $\{\theta_\alpha\}$ be
fineness endomorphisms for $(\s{C}^\bullet, d)$ with respect to $\{V_\alpha\}$.  Defining
it to be $0$ outside $\bar{V}_\alpha$ makes $\theta_\alpha s_\alpha$ a global section. 
Put $s = \sum \theta_\alpha s_\alpha$.  This is meaningful since $\theta_\alpha (x) \ne
0$ for only finitely many $\alpha$.  As $x$ has a neighbourhood meeting only finitely
many of the $U_\alpha$, it has a neighbourhood in which $s$ is a finite sum of continuous
functions and so is continuous.  But
$$d s(x) = d \sum \theta_\alpha s_\alpha(x) = \sum \theta_\alpha d s_\alpha(x) = \sum
\theta_\alpha \gamma(x) = \gamma(x).$$
We thus have $s \in \Gamma \s{C}^{n-1}$ with $ds = \gamma$. \end{proof}

\begin{lemma}\label{lem-Nfine}
If $\s{F}$ is fine then so is $(\s{N}_{\fg}^{p \bullet},\partial'')$.
\end{lemma}

\begin{proof}
Since $F$ is compact there is a cofinal set of coverings of $E$ of the
form $\varphi_\alpha(U_\alpha\times V_{\alpha i})$ where
$1\le i\le n_\alpha <\infty$ and the
$\varphi_\alpha\co U_\alpha \times F \to\pi^{-1}(U_\alpha)$ are the
coordinate functions.  We can assume that the coverings $\{U_\alpha\}$
of $B$ are locally finite. Shrink this covering to $\{U'_\alpha\}$.
For each $b \in B$, choose an index $\eta(b)$  such
that $b\in U'_{\eta(b)}$ and define endomorphisms $\theta_\alpha$
by setting for the level cochain $c$ and level $pq$-simplex $\sigma = (u, w)\co (\Delta^p,
\Delta^p \times \Delta^q) \to (B, E)$,
$$\theta_\alpha(c)(\sigma) =  \begin{cases} c(\sigma) &\text{if $\eta(u(0))=\alpha$}\\
0 &\text{otherwise,} \end{cases}$$
where $u(0)$ is the first vertex of $u$.
Then $\theta_\alpha(c)$ is $0$ outside $\pi^{-1}(\bar{U}'_\alpha)$
so $\theta_\alpha$ induces an endomorphism of $\s{N}_{\fg}^{pq}$ with
support in $\pi^{-1}(\bar{U}'_\alpha)$ and $\sum\theta_\alpha=1$.
Also $\theta_\alpha$ commutes with $\partial''$ since only the $u$ part of $\sigma$
is used in its definition.

As $\s{F}$ is fine, there are endomorphisms $\epsilon_{\alpha i}$
with support in $\bar{V}_{\alpha i}$ and $\sum_i\epsilon_{\alpha i}=1$.
Now $\s{A}|\varphi(U_\alpha\times F) \approx U_\alpha\times\s{F}$.
Let $\epsilon'_{\alpha i}$ be the endomorphism of the left hand side
corresponding to $1\times\epsilon_{\alpha i}$ on the right. Let
$\lambda_{\alpha i}$ be the endomorphism of
$\s{N}_{\fg}^{pq}=\s{N}_{\fg}^{pq}(\R)\otimes\s{A}$
defined by $\theta_\alpha\otimes\epsilon'_{\alpha i}$ on
$\pi^{-1}(U_\alpha)$ and $0$ outside of $\pi^{-1}(\bar{U}'_\alpha)$.
These patch since $\theta_\alpha=0$ outside $\pi^{-1}(\bar{U}'_\alpha)$
and clearly $\lambda_{\alpha i}$ has support in
$\varphi(U_\alpha\times V_{\alpha i})$. Also
$\sum_i\lambda_{\alpha i}=\theta_\alpha\otimes 1$, both sides being $0$
outside of $\pi^{-1}(\bar{U}'_\alpha)$ so
$\sum_{\alpha i}\lambda_{\alpha i}=1$ as required. Moreover
$\lambda_{\alpha i}$ commutes with $\partial''$.
\end{proof}

\begin{lemma}\label{lem-sNzero}
$E_1^{pq}(\s{N}_{\fg}^{pq})=0$ for $q>0$.
\end{lemma}
\begin{proof} This is local so we can assume that $E=B\times F$.   Let $\gamma \in
(\s{N}_{\fg}^{pq})_e$ and suppose $\partial'' \gamma = 0$.  Then $\gamma$ is the germ at
$e = (e_1, e_2)$ of some level $pq$-cochain $c\co S_p(U) \times S_q(V) \to \s{A}(W)$ for
some neighbourhood $W = U\times V$ of $e$.   The values $c(u,v)$ all lie in some finitely
generated submodule $A$ of $\s{A}(W)$ and we have $\partial'' c = 0$.  For every
$p$-simplex $u$ of $U$, we have a $q$-cocycle $c_u \in S^q(V; A)$.  Since $F$ is HLC,
there exists a neighbourhood $V'$ of $e_2 \in F$ such that the inclusion $V' \to V$ induces
the zero map $H^q(V; A) \to H^q(V'; A)$.  Thus we have a $(q-1)$-cochain $c'_u \in
S^{q-1}(V'; A)$ such that $\partial c'_u = c_u | V'$.  Putting $c'(u,v) = c'_u(v)$  defines a 
level $(p, q-1)$-cochain of $W' = U \times V'$ with values in $A$ and with $\partial'' c' = c|
W'$.  Passing to germs (including taking germs of elements of $A$) gives an element $\tau
\in \s{N}_{\fg}^{p,q-1}$ with $\partial'' \tau = \sigma$.
\end{proof}

\begin{lemma} Suppose the $G$-sheaf $\s{Q}$ over $F$ is fine.  Then
$E^{pq}_1(N_{\fg}(\s{Q})) = 0$ for $q > 0$. \end{lemma}
\begin{proof} This follows immediately from Lemmas \ref{lem-Nfine} and
\ref{lem-sNzero}. \end{proof}

We have now established the conditions  for dimension shifting.  We have  natural
transformations ${\psi'}^{p0}_1\co E^{p0}_1(N_{\fg}) \to E^{p0}_1(\Ler')$, $\eta^{p0}_1\co
E^{p0}_1(N_{\fg}) \to E^{p0}_1(\Ler'')$ and $i^{p0}_1\co E^{p0}_1(\Ler'') \to
E^{p0}_1(\Ler')$ induced by the inclusion $i\co \Ler'' \to \Ler'$.  These satisfy
${\psi'}^{p0}_1 = i^{p0}_1 \eta^{p0}_1$ and commute with the maps given by the first
differentials of the double complexes.  Since, for fixed $p$, $E^{pq}_1(N_{\fg})$,
$E^{pq}_1(\Ler')$ and $E^{pq}_1(\Ler'')$ are connected sequences of functors which vanish in
positive dimensions on fine sheaves over $F$, we have  natural transformations
${\psi'}^{pq}_1\co E^{pq}_1(N_{\fg}) \to E^{pq}_1(\Ler')$, $\eta^{pq}_1\co E^{pq}_1(N_{\fg})
\to E^{pq}_1(\Ler'')$ and $i^{pq}_1\co E^{pq}_1(\Ler'') \to E^{pq}_1(\Ler')$ which extend
these transformations to positive $q$.  These transformations are unique,  they commute with the
first differentials and satisfy ${\psi'}^{pq}_1 = i^{pq}_1 \eta^{pq}_1$.  Further, since
$\eta^{p0}_1$ is an isomorphism, so are the $\eta^{pq}_1$ for all $q$.

We now pass to page 2 of the spectral sequences by taking cohomology with respect to the
first differentials.  The map ${\psi'}^{pq}_2\co E^{pq}_2(N_{\fg}) \to E^{pq}_2(\Ler')$ 
induced by the map $\psi'$ of double complexes is the composite ${\psi'}^{pq}_2 = i^{pq}_2
\eta^{pq}_2$.  But $i^{pq}_2$ and $\eta^{pq}_2$ are isomorphisms.  Thus
$\psi'\co N_{\fg} \to \Ler'$ induces an isomorphism of the spectral sequences from page 2. 
Combining this with Theorem \ref{th-level} gives 

\begin{theorem}\label{th-SerL}  Suppose that the  ring $R$ is noetherian,  the
bundle $\pi\co E \to B$ has a connection,  $E$ is paracompact and HLC and that the fibre $F$ is
compact.  Let $\s{A}$ be a bundle sheaf over $E$.  Then the map $\psi' \rho \co \Ser_{\fg} \to
\Ler'$ induces an isomorphism of the Leray-Serre and Leray spectral sequences from page 2
onwards.
\end{theorem}

For A-S cohomology, the requirement that $E$ be HLC can be weakened to $B$ locally path  
connected.

\section{Comparison of Swan and Leray} \label{sec-SwL}

Let $G$ be a discrete group.  For any $G$-module $A$,  we can form the locally constant sheaf $A
\times_G EG \to BG$.  This has stalks (fibres) $A$ and the action of an element $g \in G$ on $A$ is
given by tracking the isomorphisms of fibres around a path in $BG$ corresponding to the
element $g$ of the fundamental group $G$ of $BG$.  The global section functor
$\Gamma_{BG}$ is equivalent to the $G$-invariants functor $I^G$.  Thus we may
calculate the right derived functors $R^p I^G(A)$ by resolving the sheaf $A
\times_G EG$, applying $\Gamma_{BG}$ and taking cohomology.  

Consider a $G$-sheaf $\s{F}$ on the compact $G$-space $X$.  Then by Lemma  \ref{lem-BGpara}
and Corollary \ref{cor-XGpara}, $BG$ and $X_G$ are paracompact.  By Lemma  \ref{lem-BGpara}, 
$BG$ is locally contractible and so is HLC.  It follows by  Theorem \ref{th-HLCprod} that, if $X$ is HLC,
then so is
$X_G$.  Thus we have the conditions for the use of cochain resolutions for the Swan and Leray
spectral sequences.

\textbf{Proof of Theorem \ref{th-SL}} \ \ For the Swan spectral sequence,  we first use the
A-S cochain resolution, obtaining $\bar{\s{S}}^q(X;
\s{F})$ and  take global sections $\Gamma_X\bar{\s{S}}^q(X; \s{F})$.  We then resolve as
above, using A-S cochains of $BG$.  This gives us for the Swan spectral sequence, the double complex
$$\Gamma_{BG} \bar{\s{S}}^p(BG;  \Gamma_X \bar{\s{S}}^q(X; \s{F})),$$
that is, $\ASLer^{pq}$ for the sheaf $\s{F}_G$ over $X_G = X \times_G EG$.   \qed

We are now in a position to prove Theorems \ref{th-equivar} and \ref{th-eacyc}.

\textbf{Proof of Theorem \ref{th-equivar}} \    The translation of the double
complexes gives a translation of the filtered targets of their spectral sequences.  By Theorem
\ref{th-SL}, the spectral sequences are isomorphic from page 2 onwards. Therefore their targets   
are also isomorphic.
\qed

\textbf{Proof of Theorem \ref{th-eacyc}}\ \  By assumption, $\s{F}$ is
$\Gamma^G_X$-acyclic, that is, $R^n\Gamma^G_X \s{F} = 0$ for $n > 0$.  By Theorem
\ref{th-equivar}, $H^n_G(X; \s{F}) = 0$ for $n > 0$. \qed


\begin{thebibliography}{00}
\bibitem{Bar} D. W. Barnes, \textit{Spectral sequence constructors in
algebra and topology,} Mem. Amer. Math. Soc. {53} No. 317, 1985.

\bibitem{Bred2} G. E. Bredon, \textit{Sheaf theory,} Second edition, Graduate Texts in
Mathematics {170}, Springer, New York, 1997.

\bibitem{CE} H.Cartan and S. Eilenberg, \textit{Homological Algebra,} Princeton University
Press, Princeton, 1956.

\bibitem{Dress} A. Dress, \textit{Zur Spektralsequenz von Faserungen,} Inventiones
math. {3} (1967), 172--178.

\bibitem{Dug} J. Dugundji, \textit{Topology,} Allyn and Bacon, Boston, 1970.

\bibitem{ES} S. Eilenberg and N. Steenrod, \textit{Foundations of Algebraic Topology,}
Princeton University Press, Princeton, 1952.

\bibitem{God} R.Godement, \textit{Topologie alg\'ebrique et th\'eorie des faisceaux,}
Herman, Paris, 1958.

\bibitem{Grot} A. Grothendieck, \textit{Sur quelques points d'alg\`ebre homologique,}
Tohoku Math. J. (2) {9} (1957), 119--221.

\bibitem{LW} A. T. Lundell and S. Weingram, \textit{The topology of CW complexes}, Van Nostrand
Reinhold,  New York, 1969.

\bibitem{McC} J. McCleary, \textit{User's guide to spectral sequences,} Publish or Perish
Press 1985.
 
\bibitem{Sik} A. S. Sikora, \textit{Torus and $\Z/p$ actions on manifolds,} Topology
{43} (2004), 725--748.

\bibitem{Sp} E. H. Spanier, \textit{Algebraic Topology,} McGraw-Hill, 1966.

\bibitem{Steen} N. Steenrod, \textit{The Topology of Fibre Bundles,} Princeton University
Press, Princeton, NJ, 1951.

\bibitem{Swf} R. G. Swan, \textit{A new method in fixed point theory,} Comment. Math.
Helv. {34} (1960), 1--16.

\bibitem{Sw} R. G. Swan, \textit{The theory of sheaves,} University of Chicago Press,
1964.

\end{thebibliography}
\end{document}